\documentclass{elsarticle}

\usepackage{graphicx}
\usepackage{xcolor}
\usepackage{amsmath,amssymb}

\def\eps{\varepsilon}
\def\d{\mathrm{d}}
\def\proof{\noindent {\bf Proof:~}}
\def\proofend{\hfill $ \Box $}

\newtheorem{theorem}{Theorem}[section]

\newtheorem{hypothesis}{Hypothesis}[section]
\newtheorem{definition}{Definition}[section]

\begin{document}

\begin{frontmatter}

\title{Nonautonomous control of stable and unstable manifolds in two-dimensional flows}


\author[label1,label2]{Sanjeeva Balasuriya}
\address[label1]{School of Mathematical Sciences, University of Adelaide, 
Adelaide SA 5005, Australia}
\address[label2]{Department of Mathematics, Connecticut College, New London CT
06320, USA}
\author[label3]{Kathrin Padberg-Gehle}
\address[label3]{Institute of Scientific Computing, Technische Universit\"{a}t 
Dresden, D-01062 Dresden, Germany}


\begin{abstract}
We outline a method for controlling the location of stable and unstable
manifolds in the following sense.  From a known location of the stable and unstable
manifolds in a steady two-dimensional flow, the primary segments of the manifolds
are to be moved to a user-specified time-varying location which is near the steady location.
We determine the nonautonomous perturbation to the vector field required to
achieve this control, and give a theoretical bound for the error in the manifolds 
resulting from applying this control.  The efficacy of the control strategy is illustrated
via a numerical example.
\end{abstract}


\end{frontmatter}

\section{Introduction}
\label{sec:intro}

The role of stable and unstable manifolds in demarcating flow barriers in unsteady flows
is well documented.  Determining their location in a given unsteady flow regime
is a problem which has attracted considerable attention, with many techniques continually
being developed and refined in order to improve accuracy and efficiency 
\cite{hallerberonvera,blazevskihaller,allshousethiffeault,froylandpadbergentropy,budisicmezic,hallervariational,shadden,dellnitzjunge,froylandcoherent,froylandchaos,mezicscience,tangential,open,unsteady,farazmandhaller,peacockdabiri,boffetta,krauskopf,froylandanalytic}.

Viewing this problem from the reverse viewpoint leads to an intriguing question: is it 
possible to 
force stable and unstable manifolds to lie along {\em user-defined, time-varying locations}? 
The time-variation here is arbitrarily specified, and {\em not} confined to the popular
time-periodic situation.   If possible, this would yield an invaluable tool in 
controlling transport in micro- and nano-fluidic devices, with innumerable applications.
This article answers this question in a specific setting: that of a nonautonomously perturbed
two-dimensional system, in which the issue is to determine the nonautonomous perturbation
which gives rise to the primary parts of the stable and unstable manifolds lying along 
prescribed one-dimensional curves at each instance in time.  The theory is couched in
terms of the perturbation being $ {\mathcal O}(\eps) $, and results ensuring that
the prescribed manifolds are achieved to leading-order in $ \eps $ are presented. 
Rigorous bounds for the errors in the manifolds are also established. 

The derived control strategy is tested on a time-aperiodic modification of the 
Taylor-Green flow \cite{radko,chandrasekhar,mixer,periodic,l2mixer}.   Numerical
diagnostics are compared with the prescribed stable manifolds, and excellent results are obtained.

 While the method developed in this article is confined to perturbations of autonomous flows, it is to our knowledge the first theoretical contribution
towards developing a control strategy for stable and unstable manifolds in nonautonomous
flows.  As such, it may serve as an important initial step towards building a more complete theory
for the nonautonomous control of flow barriers.

\section{Controlling stable manifold}
\label{sec:stablecontrol}

Consider for $ x \in \Omega $, a two-dimensional open connected set, the system
\begin{equation}
\dot{x} = f(x)
\label{eq:unp}
\end{equation}
in which $ f: \Omega \rightarrow {\mathbb R}^2 $, and sufficient smoothness will be
assumed (to be characterised shortly).

\begin{figure}[t]
\includegraphics[scale=1.0]{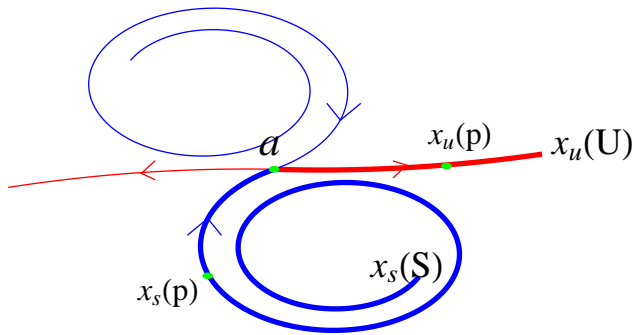}
\caption{The restricted stable manifold $ \Gamma_s $ (thick blue curve) and unstable 
manifold $ \Gamma_u $ (thick red curve) of the fixed point $ a $, in the situation in which
the stable manifold exhibits spiralling
and the unstable manifold escapes to infinity.}
\label{fig:unp}
\end{figure}

\begin{hypothesis}[Saddle point at $ a $]
The system (\ref{eq:unp}) possesses a saddle fixed point $ a $, that is, $ f(a) = 0 $ and
$ D f(a) $ possesses real eigenvalues $ \lambda_s $ and $ \lambda_u $ such that
$ \lambda_s < 0 < \lambda_u $.  
\label{hyp:saddle}
\end{hypothesis}
Then, $ a $ possesses corresponding one-dimensional
stable and unstable
manifolds.  We will focus on {\em segments} of {\em one branch} of each of these manifolds, and denote
them  by $ \Gamma_s $ and $ \Gamma_u $ respectively.  The
segment of the stable manifold branch we will consider can be represented
parametrically by 
\[
\Gamma_s := \left\{ x_s(p) ~: ~ p \in [S, \infty) \right\} 
\]
in which $ x_s(t) $ is a solution to (\ref{eq:unp}) with initial condition $ x_s(0) 
\in \Gamma_s $, and $ S \in \left( - \infty, 0 \right] $ represents a  finite backwards 
time until which the trajectory is evolved.
Notice in particular that $ x_u(t) \rightarrow a $ as $ t \rightarrow  \infty $, and so
$ \bar{\Gamma}_s $ contains $ a $, while the other end of the curve segment comprising
$ \Gamma_s $ ends at the point $ x_s(S) $.  From this definition, it is clear that
$ \Gamma_s $ {\em cannot} be (i) a branch of a stable manifold which has infinite length, or
(ii) a heteroclinic or homoclinic manifold associated with a fixed point since $ x_s(t) $ cannot approach a fixed point in finite time.
On the other hand, $ \Gamma_s $ could be any other finite length restriction of a branch of the
stable manifold emanating from $ a $, including a segment of any of the above two entities,
or a segment of a manifold which has many rotations as it spirals out from a limit cycle. 
  Similarly,
let $ \Gamma_u $ be a restricted branch of the unstable manifold of $ a $ which is
parametrisable as
\[
\Gamma_u := \left\{ x_u(p) ~: ~ p \in (-\infty, U] \right\} 
\]
in which $ x_u(t) $ is a solution to (\ref{eq:unp}) with initial condition
$ x_u(0) \in \Gamma_u $, and which satisfies $ x_u(t) \rightarrow a $ as 
$ t \rightarrow - \infty $, and $ U \in [0, \infty) $ is a finite forward time until which the
trajectory is evolved.  See Fig.~\ref{fig:unp} for an example of the finite segments $ \Gamma_s $ and
$ \Gamma_u $.

The goal is to determine a nonautonomous perturbation to the vector field in the form
\begin{equation}
\dot{x} = f(x) + \eps g(x,t)
\label{eq:pert}
\end{equation}
in which $ \eps \in [0, \eps_0) $ where $ \eps_0 \ll 1 $, 
such that 
$ \Gamma_s $ and $ \Gamma_u $ perturb to $ \eps $-close time-dependent entities
which are {\em specified}.  
The following smoothness hypotheses on the functions
$ f $ and $ g $ will be assumed, in which $ D $ represents the spatial (matrix) derivative
operator in $ \Omega $.

\begin{hypothesis}[Smoothness of $ f $ and $ g $]
The functions $ f : \Omega \rightarrow {\mathbb R}^2 $ and $ g: \Omega \times
{\mathbb R} \rightarrow {\mathbb R}^2 $ satisfy the following smoothness and 
boundedness assumptions.
\begin{itemize}
\item[(f)] $ f \in {\mathrm{C}}^2 \left( \Omega \right) $, and is 
such that there exists a constant
$ C_f $ satisfying
\begin{equation}
\left\| f \left( x \right) \right\| + \left\| D f \left( x \right) \right\| + \left\| D^2 f \left( x \right)
\right\| \le C_f ~~{\mathrm{for~all}}~x \in \Omega \, .
\label{eq:fbound}
\end{equation}
\item[(g)] $ g \in {\mathrm{C}}^2 \left( \Omega \right) $ for each $ t \in {\mathbb R} $, and
$ g \in {\mathrm{C}}^1 \left( {\mathbb R} \right) $ for each $ x \in \Omega $, and moreover
there exists a constant $ C_g $ satisfying
\begin{equation}
\left\| g \left( x,t  \right) \right\| + \left\| D g \left( x, t \right) \right\| + \left\|\frac{\partial g}{\partial t} \left( x,t  \right)
\right\| \le C_g ~~{\mathrm{for~all}}~(x,t) \in \Omega 
\times {\mathbb R} \, .
\label{eq:gbound}
\end{equation}
\end{itemize}
\label{hyp:fgbound}
\end{hypothesis}
A note on the norms used in (\ref{eq:fbound}) and (\ref{eq:gbound}) is in order.  The
norm on $ {\mathbb R}^2 $ is the standard Euclidean norm; we could have stated the
relevant norms on $ f $, $ g $ and $ \partial g/\partial t $ by using the modulus $ \left| \centerdot \right| $ instead.  The norm on
the $ 2 \times 2 $ matrices $ D f $ and $ D g $ is the operator norm induced by the
Euclidean norm.  The norm on the $ 2 \times 2 \times 2 $ entity $  D^2 f  $ is the induced
operator norm associated with the above norms on vectors and matrices, i.e., 
\begin{equation}
\left\| \left(D^2 f  \right) \right\| = \sup_{v \in {\mathbb R}^2 \setminus 0} \frac{\left\| 
\left( D^2 f \right) \, v \right\|}{\left\| v \right\|} \, , 
\label{eq:operatornorm}
\end{equation}
in which the previously mentioned operator norm definition for $ 2 \times 2 $ matrices 
is used in the numerator.

Now, the smooth function $ g(x,t) $ will be the control which achieves the
desired stable and unstable manifolds, which can be now represented by
$ \Gamma_s^\eps(t) $ and $ \Gamma_u^\eps(t) $ respectively.  In viewing these
restrictions to the manifolds in this nonautonomous setting, it makes sense to 
represent (\ref{eq:pert}) in the augmented form
\begin{equation}
\left. \begin{array}{l}
\dot{x} = f(x) + \eps g(x,t) \\
\dot{t} = 1
\end{array} \right\} \, 
\label{eq:pertaug}
\end{equation}
with phase space now being $ \Omega \times {\mathbb R} $.
For (\ref{eq:pertaug}) when $ \eps = 0 $, the conditions stated for (\ref{eq:unp}) provide
for the presence of a hyperbolic trajectory $ \left( a, t \right) $ with two-dimensional
stable and unstable manifolds.  

\begin{figure}[t]
\includegraphics[height=7cm, width=12cm]{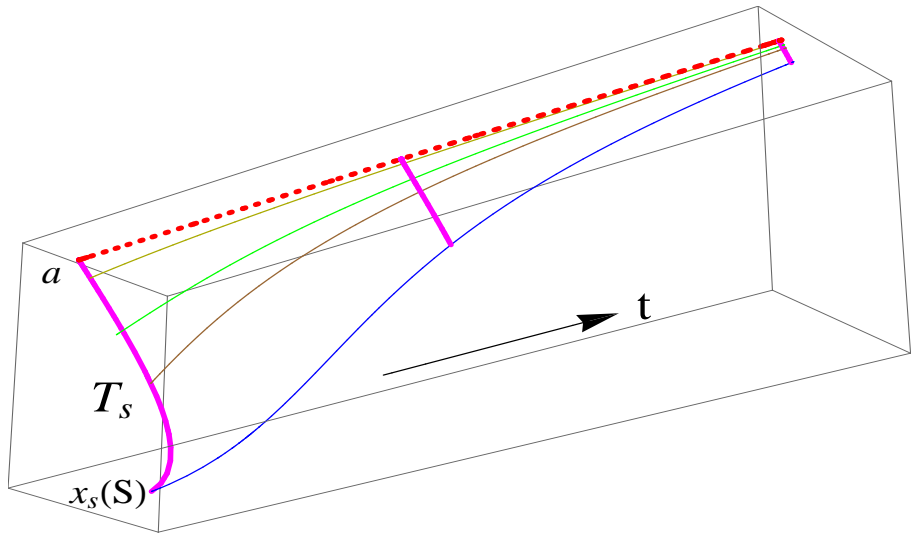}
\caption{The evolution of the restricted stable manifold curve (heavy magenta
curves), and the hyperbolic trajectory (dashed red line).}
\label{fig:shortening}
\end{figure}

From this point onwards, this Section will focus only on controlling the {\em stable} manifold,
with the unstable manifold control description postponed to the subsequent Section.
It will be necessary to restrict the stable manifold in 
time in the following sense.  Let $ T_s < 0 $ be a time-value beyond which the restricted stable
manifold is to be defined.  This signs for $ T_s $ is chosen to
ensure that $ t = 0 $ is a legitimate choice for both the restricted stable and
unstable manifolds.  Restricting time in this way will be necessary when $ \eps \ne 0 $ because the restrictions on $ p $ mean that only segments of the relevant manifolds
are defined in each time-slice, and since these segments evolve with time further
restrictions will arise. 
The restricted two-dimensional
stable manifold of (\ref{eq:pertaug}) when $ \eps = 0 $ will be
represented in parametric form by
\begin{equation}
\Gamma_s := \left\{ \left( x_s(t-T_s+p),t \right) ~: ~ (p,t) \in [S, \infty) \times 
[T_s, \infty)  \right\} \, , 
\label{eq:unpmanifolds} 
\end{equation}
in which the notation $ \Gamma_s $ is retained with an abuse of notation.  This
parametrisation with respect to $ (p,t) $ has been chosen such that the parameter
$ p $ selects the specific trajectory on the relevant manifold of (\ref{eq:pertaug}) when
$ \eps = 0 $, and $ t $ represents the time-evolution of that trajectory.  
So for example
if a point $ x_s(p) $ is chosen in the time-slice $ t = T_s $, then $ \left( x_s(t-T_s+p), t \right) $
represents the corresponding forward trajectory on $ \Gamma_s $ as it evolves with time
$ t $.
In the time-slice
$ T_s $, the restriction $ p \ge S $ implies
the relevant segment of the unperturbed stable manifold goes from $ x_s(S) $
to $ a $.  Now,
in a time-slice $ t > T_s $, the trajectory through $ x_s(S) $ would have evolved to
the location $ x_s(S+t-T_s) $, and information on $ x_s(p) $ for $ p $ values
less than $ S+t-T_s $ cannot be available since such would correspond to points on
the stable manifold which were beyond $ x_s(S) $ at time $ T_s $.  As time $ t $ evolves
for the unperturbed steady flow (\ref{eq:unp}),
$ x_s(S+t-T_s) $ approaches the saddle fixed point $ a $, and therefore the length of
the restricted stable manifold in each time-slice gets shorter.  In other words, the
restrictions implied in (\ref{eq:unpmanifolds}) results are actually associated with 
shorter and shorter segments of the stable manifold as time gets larger.
This is illustrated in Fig.~\ref{fig:shortening}, in which trajectories associated with five
$ p $ values are shown beginning with an ``initial'' point in the time-slice $ T_s $.
The ``furthest'' of these corresponds to the initial point $ x_s(S) $ (i.e., $ p = S $),
and after three intermediate $ p $ values, the dashed trajectory 
is $ (a,t) $, which can be thought of as  $ p = \infty $ in (\ref{eq:unpmanifolds}).
The stable manifold in the time-slice $ T_s $ is the curve segment (heavy magenta curve)
which connects
together all five starting points.   In time-slices $ t $ as time evolves, all trajectories
get closer together (indeed, they get closer to the dashed hyperbolic trajectory 
$ (a,t) $), which means that the restricted stable manifold are becoming curves
of
smaller and smaller length in each time-slice. These are indicated at two later time values
in Fig.~\ref{fig:shortening}, also as heavy magenta curves. A similar description 
works for
the restricted unstable manifolds in (\ref{eq:unpmanifolds}), but in this case the
manifold segments in each time-slice becomes shorter curves as $ t $ decreases.

Now, when $ \eps \ne 0 $, and for any $ g $ satisfying the smoothness assumptions
in Hypothesis~\ref{hyp:fgbound}, 
the hyperbolic trajectory $ \left( a, t \right) $ of (\ref{eq:pertaug})
perturbs to an $ {\mathcal O}(\eps) $-close trajectory $ \left( a_\eps(t), t \right) $ which retains hyperbolicity.
The proof of this is via exponential dichotomies \cite{coppel,yi,yagasakicontrol}, and as a consequence this trajectory
retains stable and unstable manifolds which are $ \eps $-close to the original ones.  In particular, it retains
a stable manifold $ \eps $-close to (\ref{eq:unpmanifolds}).  The locations of the manifold will depend on
the choice of $ g $, but here we {\em specify} the perturbed restricted manifolds, and
find additional conditions on $ g $ in order to achieve these.  The desired restricted
stable manifold will be
represented parametrically by
\begin{equation}
\Gamma_s^\eps := \left\{ \left( x_s^\eps (p,t),t \right) ~:~ (p,t) \in [S,\infty) \times [T_s,\infty)  \right\}  \, , 
\label{eq:pertmanifolds} 
\end{equation}
where  the parametrisation
$ x_s^\eps $ is assumed {\em given}, but satisfies several conditions to ensure
consistency.  To express these conditions,  we first define
\begin{equation}
J = \left( \begin{array}{lr}
0 & - 1 \\ 1 & 0 
\end{array} \right) \, ,
\label{eq:J}
\end{equation}
the premultiplicative matrix which rotates vectors in $ {\mathbb R}^2 $ by $ + \pi / 2 $.

\begin{hypothesis}[Stable manifold requirements]
For each $ t \ge T_s $,  the quantity 
$ \left\{ x_s^\eps(p,t)~:~p \ge S \right\} $ is a curve in 
$ \Omega $.  These restricted stable manifold curves satisfy the following conditions.
\begin{itemize}
\item[(a)]  [Smoothness]  There exists a constant $ K_s > 0 $ such that for all
$ (p,t) \in [S,\infty) \times [T_s,\infty) $ and for all $ \eps \in (0,\eps_0) $,
\begin{small}
\begin{equation}
 \hspace*{-1cm} \left|  x_s^\eps(p,t) \right| + \left| \frac{\partial}{\partial t} x_s^\eps(p,t) \right| +
\left|  \frac{\partial}{\partial p} x_s^\eps(p,t) \right| + \left| \frac{\partial}{\partial \eps} x_s^\eps(p,t) \right| + \left| \frac{\partial^2}{\partial \eps^2} x_s^\eps(p,t) \right|
+ \left| \frac{\partial^3}{\partial \eps^3} x_s^\eps(p,t) \right| < K_s \, .
\label{eq:xssmooth} 
\end{equation}
\end{small}
\item[(b)] [Closeness]  There exists a constant $ C_s > 0 $ such that for all
$ (p,t) \in [S,\infty) \times [T_s,\infty) $, 
\begin{equation}
\left| x_s^\eps (p,t) - x_s(t-T_s+p) \right| + \left| 
\frac{\partial}{\partial p} \left(x_s^\eps (p,t) - x_s(t-T_s+p) \right) \right|
\le C_s \, \eps \, . 
\label{eq:xsclose} 
\end{equation}
\item[(c)] [Limit]  For each $ t \ge T_s $,  $ 
\displaystyle \lim_{p \rightarrow \infty} x_s^\eps(p,t) $ is well defined.
\item[(d)] [Mappability] For each $ t \ge T_s $, there exist intervals
$ \left[ S_1(t), S_2(t) \right) $ and $ \left[ S_1^\eps(t), S_2^\eps(t) \right) $ ---both
of which are contained in $ [S, \infty) $---  and 
 a scalar function $ r_s(\centerdot,t) $ defined on $ \left[ S_1(t), S_2(t) \right) $ 
 which satisfies 
 \begin{equation}
 x_s^\eps(p^\eps,t) = x_s(t-T_s+p) + r_s(p,t)  \frac{J f \left( x_s(t-T_s+p) \right) }{ \left|  J f \left( x_s(t-T_s+p) \right)  \right|} \, ,
 \label{eq:mappables}
 \end{equation}
 such that the mapping $ p \rightarrow p^\eps $ 
 from $ \left[ S_1(t), S_2(t) \right) $ to $ \left[ S_1^\eps(t), S_2^\eps(t) \right) $ 
 defined through (\ref{eq:mappables})
 is a diffeomorphism.
\item[(e)] [Congruence at time zero]   The $ p $ parameters in the time-slice $ t = 0 $ between
the unperturbed and the required restricted stable manifold curves match up, i.e., for
all $ p \in \left[ S_1(0), S_2(0) \right) $,  
\begin{equation}
\left[f \left( x_s(-T_s+p) \right) \right]^T  \left[ x_s^\eps(p,0) - x_s(-T_s+p) \right] = 0 \, .
\label{eq:xscongruence}
\end{equation}
\end{itemize}
\label{hyp:xs}
\end{hypothesis}

\begin{figure}[t]
\includegraphics[scale=1.2]{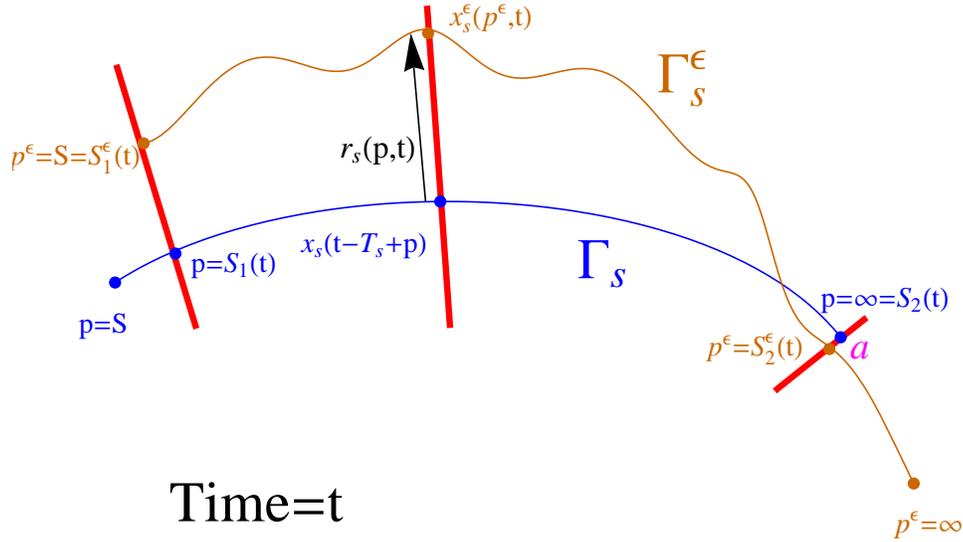}
\caption{An illustration of the mappability condition (\ref{eq:mappables}) in the time-slice $ t $.  The heavy lines are in the normal direction $ J f $ to $ \Gamma_s $.  
 The interval $ \left[ S_1(t), S_2(t) \right) $ is the $ p $-interval for which the it is 
possible to map from $ \Gamma_s $ to $ \Gamma_s^\eps $ by going in the normal 
direction $ \Gamma_s $, while $ \left[ S_1^\eps(t), S_2^\eps(t) \right) $ is the
corresponding interval for $ p^\eps $ which parametrises $ \Gamma_s^\eps $.  In this
pictured situation, $ S_2(t) = \infty $ and $ S_1^\eps(t) = S $.}
\label{fig:mappables}
\end{figure}

These hypotheses require some explanation.  The condition (\ref{eq:xssmooth}) is
a straightforward requirements on the smoothness and boundedness
of the required restricted manifold.  The condition (\ref{eq:xsclose}) is a
$ {\mathcal O}(\eps) $-closeness requirements between 
$ x_s^\eps(p,t) $ and $ x_s(t-T_s+p) $ at each $ (p,t) $ value.  In particular,
for each fixed $ t \in [T_s,\infty) $, the curves 
$ x_s^\eps(p,t) $ and $ x_s(t-T_s+p) $ and their tangents in the $ t $ time-slice 
are assumed to remain $ \eps $-close.  Condition~(c) requires that the
end of the curve -- that purportedly is on the hyperbolic trajectory $ a_\eps(t) $ -- 
is well-defined.  While becoming unbounded is already precluded by condition~(a),
condition~(c) prevents $ x_s^\eps(p,t) $ behaving like, say, $ \cos p $ for large $ p $.  

The condition in Hypothesis~\ref{hyp:xs}(d) prevents for example choosing $ x_s^\eps(p,t) $ such that
a self-intersecting curve is generated in a time-slice.  The intuition is that in each time
slice $ t $ the restricted
autonomous stable manifold segment
and the required restricted nonautonomous stable manifold segment 
are mappable to one another by proceeding in the normal direction to each point $ x_s(t-T_s-p) $, by a signed distance $ r_s(p,t) $.  The restriction of $ p $ and $ p^\eps $ 
to these subintervals of $ [S,\infty) $ is since some parts of the required $ \Gamma_s^\eps $
may venture ``beyond'' the span of the normal direction to $ \left\{ x_s(t-T_s+p) ~;~
p \in [S, \infty) \right\} $.  This condition is illustrated by example in Fig.~\ref{fig:mappables}.
This mapping from $ \Gamma_s $ to $ \Gamma_s^\eps $ by going along the normal
direction $ J f $ from each point on $ \Gamma_s $ parametrised by $ p $ must be a
diffeomorphism from $ \left[ S_1(t), S_2(t) \right) $ to $ \left[ S_1^\eps(t), S_2^\eps(t) 
\right) $, which prevents $ \Gamma_s^\eps $ having self-intersections or twists which
make the inverse function undefined.

\begin{figure}[t]
\includegraphics[scale=1.2]{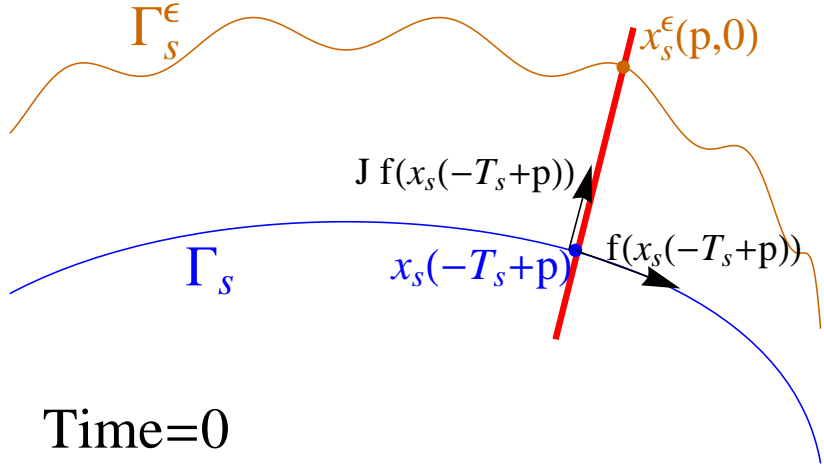}
\caption{The congruence condition (\ref{eq:xscongruence}) in the time-slice $ t = 0 $: 
at each $ p $, the normal vector at $ x_s(-T_s+p) $ meets $ \Gamma_s^\eps $ 
at the point $ x_s^\eps(p,0) $.}
\label{fig:congruence}
\end{figure}

Finally, the congruence condition~(e) reflects a choice of parametrisation taken 
in the time-slice $ t = 0 $, which is shown in Fig.~\ref{fig:congruence} for
the restricted stable manifold.   For any fixed $ p $, consider the point $ x_s(-T_s+p) $
on the unperturbed stable manifold, and suppose we draw a line perpendicular 
to $ f \left( x_s(-T_s+p) \right) $ in this time-slice $ t = 0 $, as shown in Fig.~\ref{fig:congruence}.  Now, the intersections of $ \Gamma_s^\eps $ and $ \tilde{\Gamma}_s^\eps $ in the time-slice $ 0 $ are also shown in this figure, and the normal line intersects each of these curves.  The congruence condition (\ref{eq:xscongruence}) in the time-slice
$ t = 0 $ means that
the $ p $-parametrisation of $ x_s^\eps(p,0) $ is chosen such that $ x_s^\eps(p,0) $
lies exactly on this normal line drawn at $ x_s(-T_s+p) $.  We have the freedom to do
this for all mappable $ p $ in this one particular time-slice; it is merely a choice of parametrisation
of the one dimensional curve obtained by intersecting $ \Gamma_s^\eps $ with the time-slice
$ \left\{ t = 0 \right\} $.

\begin{figure}[t]
\includegraphics[scale=1.2]{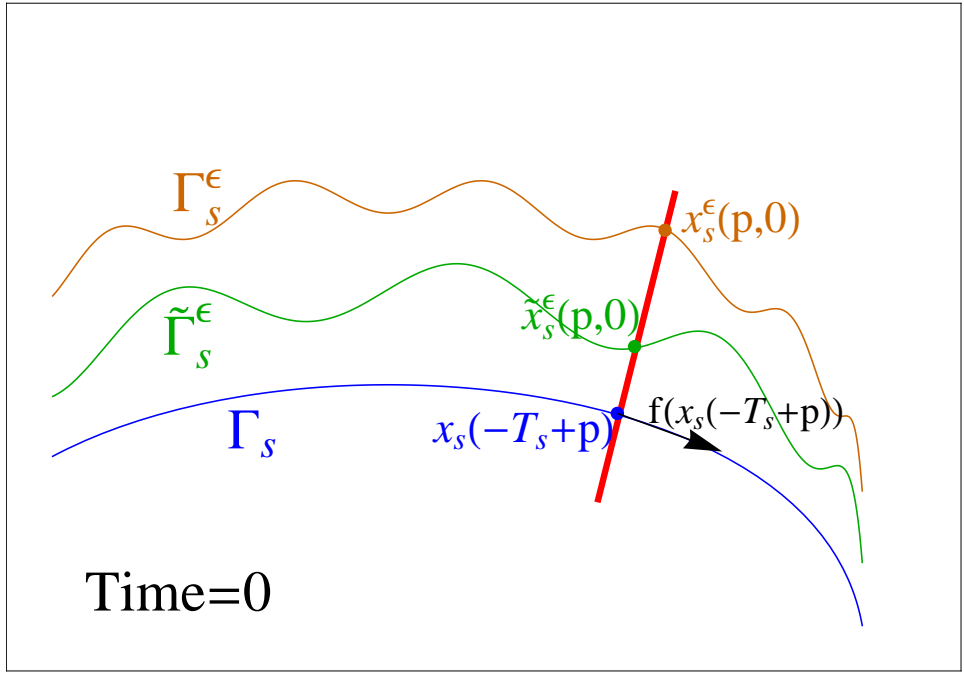}
\caption{The congruence conditions (\ref{eq:xscongruence}) and (\ref{eq:xtildescongruence})
in the time-slice $ t = 0 $, illustrating that both the required ($\Gamma_s^\eps $)
and the real ($\tilde{\Gamma}_s^\eps $) stable manifolds have a congruent
$ p $-paramatrisation at time zero.}
\label{fig:congruence+real}
\end{figure}

While the desired restricted stable manifold is given by  (\ref{eq:pertmanifolds}), 
Hypotheses~\ref{hyp:xs} further restricts the $ (p,t) $ values to lie in the set
\begin{equation}
\Xi_s := \left\{ \left( p,t \right) ~:~ t \ge T_s~{\mathrm{and}}~ S_1(t) \le p < S_2(t)  \right\}  \, . 
\label{eq:ptsets} 
\end{equation}
We will assume that the largest interval $ \left[ S_1(t), S_2(t) \right) $ 
has been chosen for each $ t $ in order to fulfil the 
mappability condition of Hypothesis~\ref{hyp:xs}.  

For $ (p,t,\eps) \in \Xi_s \times (0,\eps_0) $, we define
\begin{equation}
M_s^\eps(p,t) := 
 \left[ J f \left( x_s(t-T_s+p) \right) \right]^T \frac{x_s^\eps (p,t) - x_s(t-T_s+p)}{\eps}
\label{eq:ms}
\end{equation}
and
\begin{equation}
B_s^\eps(p,t) := \left[ f \left( x_s(t-T_s+p) \right) \right]^T \frac{x_s^\eps (p,t) - x_s(t-T_s+p)} {\eps} \, ,
\label{eq:bs}
\end{equation}
which respectively represent projections of the difference between the unperturbed and
the desired restricted stable manifold in the normal and tangential directions to the original 
manifold $ \Gamma_s $ in the time-slice $ t $ .
Note that for a specified $ x_s^\eps(p,t) $, both $ M_s^\eps $ and $ B_s^\eps $ can be
computed numerically based on the above expressions.
Now, the required values of $ g $ (to leading-order) 
shall be expressed in terms of an orthogonal
basis formed by projecting normally and tangentially to the autonomous stable
manifold at $ x_s(t-T_s+p) $ in the time-slice $ t $.

\begin{definition}[Control velocity for stable manifold]
The control velocity $ g $ satisfies the smoothness conditions of Hypothesis~\ref{hyp:fgbound}, and 
moreover is specified by
\begin{eqnarray}
g^\perp \left( x_s(t-T_s+p),t \right)  & := & \frac{\left[ J f \left( x_s(t-T_s+p) \right) \right]^T}{
\left|  f \left( x_s(t-T_s+p) \right) \right|} g \left( x_s(t-T_s+p), t \right) \nonumber \\
& = & \frac{\frac{\partial M_s^\eps}{\partial t} (p,t) - {\mathrm{Tr}} \, 
\left( D f  \right) M_s^\eps(p,t)}{\left| f  \right|}
\label{eq:gsperp}
\end{eqnarray}
and
\begin{eqnarray}
g^{\|} \left( x_s(t-T_s+p),t \right) & := & \frac{\left[ f \left( x_s(t-T_s+p) \right) \right]^T}{
\left|  f \left( x_s(t-T_s+p) \right) \right|} g \left( x_s(t-T_s+p), t \right) \nonumber \\
& &  \hspace*{-2.5cm} = \frac{ \left| f \right|^2 \frac{\partial B_s^\eps(p,t)}{\partial t}
- f^T \left[ (D f) + (D f)^T \right] \left[ J f M_s^\eps(p,t) + f B_s^\eps(p,t)\right]}{ \left| f \right|^3} 
\label{eq:gsparallel} 
\end{eqnarray}
in which $ f $ and $ D f $ in the above expressions are evaluated at $ x_s(t-T_s+p) $,
\label{def:stable}
\end{definition}

By choosing $ g $ as above, it will be possible to achieve the desired nonautonomous stable
manifold correct to $ {\mathcal O}(\eps) $.  We will in Theorem~\ref{theorem:stable}
specify the error precisely.   First, let us describe how to apply this control velocity computationally to
achieve the desired stable manifold. Suppose we are given
the parametrised form $ x_s^\eps(p,t) $ of the restricted manifold $ \Gamma_s^\eps $,
and full knowledge of the nearby unperturbed steady flow (\ref{eq:unp}).  To compute
the control condition required to obtain the restricted stable manifold to leading-order,
we proceed as follows.
\begin{enumerate}
\item Since full knowledge of the unperturbed steady flow (\ref{eq:unp}) is presumed known,
compute $ x_s(p) $, and hence compute $ f \left( x_s(t-T_s+p) \right) $, $ D f \left( x_s(t-T_s+p) \right) $ and $
{\mathrm{Tr}} \, D f \left( x_s(t-T_s+p) \right) $  as functions of $ (p,t) $;
\item Since the restricted perturbed manifold $ \Gamma_s^\eps $ is presumed
specified through its parametrisation $ x_s^\eps(p,t) $, 
compute $ M_s^\eps(p,t) $ and $ B_s^\eps(p,t) $ from (\ref{eq:ms}) and 
(\ref{eq:bs}), recalling the restriction $ (p,t) \in \Xi_s$;
\item Determine the $ t $-derivatives of both $ M_s^\eps(p,t) $ and $ B_s^\eps(p,t) $,
using a numerical method if needed;
\item Substitute these values into (\ref{eq:gsperp}) and (\ref{eq:gsparallel})
to determine $ g_s^\perp $ and $ g_s^\| $, where
in each time-slice $ t $, the values are found along the
restricted part of $ \Gamma_s $ lying between $ x_s(t-T_s+S) $ and $ a $;
\item Since $ g^\perp $ and $ g^\| $ give the components of $ g $ in the directions $
J f $ and $ f $ respectively, this determines $ g $ at the locations $ x_s(t-T_s+p) $ in time-slices $ t $;
\item Extend $ g $ in any suitably relevant fashion to the spatial domain while being
consistent with this requirement.
\end{enumerate}

To characterise the fact that this procedure results in a nonautonomous stable manifold which is correct
to $ {\mathcal O}(\eps) $, and to additionally quantify the error resulting from this process, we need to compare the
{\em desired}
stable manifold as specified in Hypothesis~\ref{hyp:xs} with the 
{\em true} stable manifold resulting from applying the control velocity of Def.~\ref{def:stable}.
We define this true stable manifold by
\begin{equation}
\tilde{\Gamma}_s^\eps := \left\{ \left( \tilde{x}_s^\eps (p,t),t \right) ~:~ (p,t) \in \Xi_s  \right\}  \, , 
\label{eq:pertmanifoldstilde} 
\end{equation}
rather than (\ref{eq:pertmanifolds}), in which the $ \tilde{x}_s^\eps(p,t) $ 
is an exact trajectory of (\ref{eq:pert}) in which $ g $ is as specified in Def.~\ref{def:stable}.  
For each $ p $, the trajectory $ \tilde{x}_s^\eps(p,t) $  lies on the associated true perturbed manifold 
$ \tilde{\Gamma}_s^\eps $, with the $ t $ parametrising the time evolution.   Thus, $ \left| \tilde{x}_s^\eps(p,t) -
a_\eps(t) \right| \rightarrow 0 $ as $ t \rightarrow \infty $ for any $ p $.
Moreover, the parametrisation $ p $ can be chosen so that $ \tilde{x}_s^\eps(p,t) $
is $ {\mathcal O}(\eps) $-close to $ x_s(t-T_s+p) $, that is, there exists 
 a constant $ \tilde{C}_s $ such that
\begin{equation}
\left| \tilde{x}_s^\eps (p,t)  - x_s(t-T_s+p) \right| +
\left| \frac{\partial}{\partial p} \left( \tilde{x}_s^\eps (p,t)  - x_s(t-T_s+p)
\right) \right| \le
\eps \tilde{C}_s
\label{eq:closenesss}
\end{equation}
for $ (p,t,\eps) \in \Xi_s \times [0,\eps_0) $. The expectation
is that $ \tilde{C}_s \approx C_s $ as given in (\ref{eq:xsclose}), since the purported
$ \Gamma_s^\eps $ as given in (\ref{eq:pertmanifolds}) and parametrised by
$ x_s^\eps(p,t) $ will be forced to be close
to the true restricted stable manifold $ \tilde{\Gamma}_s^\eps $ 
which is parametrised by $ \tilde{x}_s^\eps(p,t) $.  
We note from Fig.~\ref{fig:congruence+real} that it is possible to choose the parametrisation
$ p $ on $ \tilde{x}_s^\eps(p,t) $ such that it too lies exactly on the normal vector drawn
at $ x_s(-T_s+p) $ in the time-slice $ t = 0 $.  Essentially, we can {\em choose} the
points $ \tilde{x}_s^\eps(p,0) $ 
on the normal vector as initial conditions for (\ref{eq:pert}), thereby  defining the parameter values $ p $ which identify each trajectory in this way.  That is, analogous to the
congruence condition (\ref{eq:xscongruence}) at time zero for the desired stable manifold, we require that 
\begin{equation}
\left[f \left( x_s(-T_s+p) \right) \right]^T  \left[ \tilde{x}_s^\eps(p,0) - x_s(-T_s+p) \right] = 0
\label{eq:xtildescongruence}
\end{equation}
for the {\em true} stable manifold.
(For more details about characterising
such tangential movement of perturbed manifolds, see \cite{tangential}.)

\begin{figure}[t]
\includegraphics[scale=1.2]{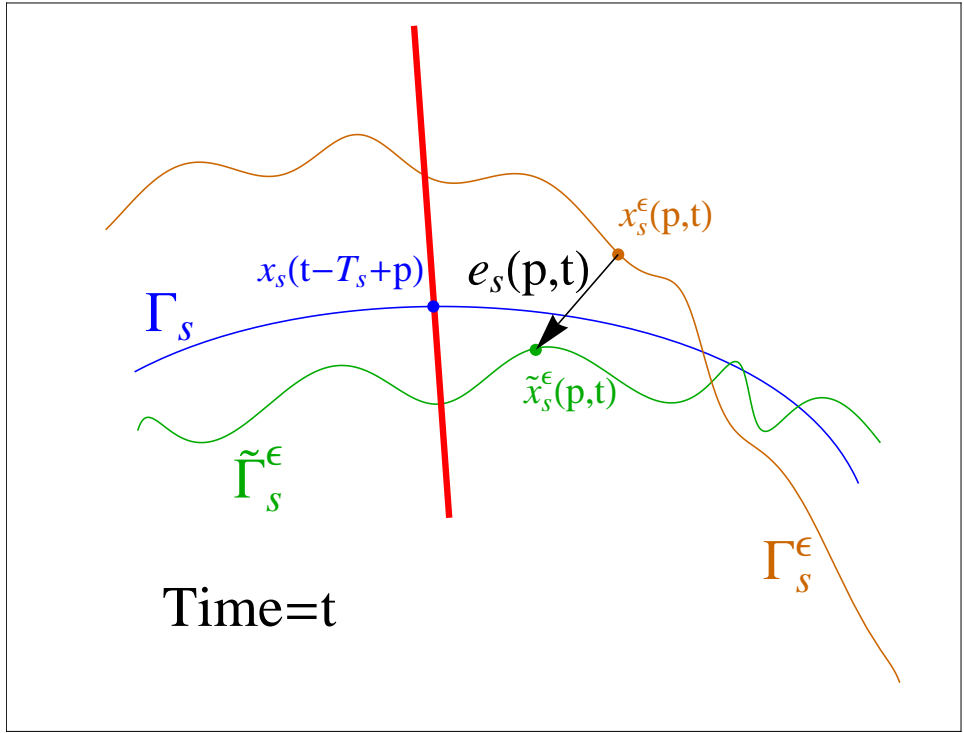}
\caption{The intersections of the unperturbed ($ \Gamma_s $), required ($\Gamma_s^\eps$)
and true ($\tilde{\Gamma}_s^\eps $) restricted stable manifolds in a general time-slice $ t $.}
\label{fig:timeslice}
\end{figure}

Now, we write
\begin{equation}
\tilde{x}_s^\eps(p,t) =  x_s^\eps(p,t) + e_s(p,t)  \, ,
\label{eq:xstilde}
\end{equation}
in which the $ e_s(p,t) $s represent the error in the restricted stable manifold at
time $ t $ and associated with the parametrisation $ p $.  An illustration of $ e_s(p,t) $
is provided in Fig.~\ref{fig:timeslice}.  We note that while in the time-slice $ t = 0 $ 
the parameter $ p $ was chosen to ensure that all three points corresponding to the
same parameter value $ p $ lie on the normal line to the unperturbed manifolds drawn
at $ x_s(t-T_s+p) $, this is not necessary so in a general time-slice.  This is because 
the $ t $-evolution of $ \tilde{x}_s^\eps(p,t) $ is generated by the flow (\ref{eq:pert}), and because the 
$ t $-evolution of $ x_s^\eps(p,t) $ is {\em specified}.  Thus, the error
term $ e_s(p,t) $ has in general both a normal and a tangential term.  Bounds for these
components of the error can be stated precisely as follows.

\begin{theorem}[Error in stable manifold]
Assume the control velocity $ g $ satisfies Def.~\ref{def:stable}, and define
\begin{equation}
e_s^\perp(p,t) := \frac{J f \left( x_s(t-T_s+p) \right)}{\left| f \left( x_s(t-T_s+p) \right) \right|}
e_s(p,t)
~~{\mathrm{and}}~~
e_s^\|(p,t) := \frac{ f \left( x_s(t-T_s+p) \right)}{\left| f \left( x_s(t-T_s+p) \right) \right|} 
e_s(p,t) \, .
\label{eq:errorstabledef}
\end{equation}
The normal component is bounded by
\begin{small}
\begin{equation}
 \left| e_s^\perp(p,t) \right| \le  \left[ C_s C_g + \frac{C_s^2 C_f}{2} \right] \eps^2 \frac{\int_t^\infty 
\left| f \left( x_s(\tau\!-\!T_s\!+\!p) \right) \right|  
\exp \left[ \int_\tau^t  {\mathrm{Tr}} \left[ D f \left( x_s(\xi\!-\!T_s\!+\!p) \right) \right] \d \xi \right] 
\d \tau}{\left| f \left( x_s(t\!-\!T_s\!+\!p) \right) \right|} \, ,
\label{eq:stableerrorboundperp}
\end{equation}
\end{small}
 for $ (p,t,\eps) \in \Xi_s \times (0,\eps_0) $, and  satisfies
 the limits
\begin{equation}
\lim_{t \rightarrow \infty} \left| e_s^\perp(p,t) \right| \le - \frac{ \left[ C_s C_g + \frac{C_s^2 C_f}{2} \right] \eps^2}{\lambda_s} \, , \, 
\lim_{p \rightarrow \infty} \left| e_s^\perp(p,t) \right| \le  \frac{ \left[ C_s C_g + \frac{C_s^2 C_f}{2} \right] \eps^2}{\lambda_u}  \, ,
\label{eq:stableerrorboundperpinfinity}
\end{equation}
as long as these limits can be taken within the domain $ \Xi_s $.
The tangential component of the error is
bounded by
\begin{small}
\begin{eqnarray}
\left| e_s^\|(p,t) \right| & \le &  \eps^2 \left( C_s C_g + \frac{ C_s^2 C_f}{2} \right) \left| f \left( x_s(t-T_s+p) \right) \right| \nonumber \\
& & \hspace*{-2.4cm} \times \left| \int_0^t 
\frac{ \left| f \left( x_s(\tau\!-\!T_s\!+\!p) \right) \right| + 2 C_f \int_\tau^\infty 
\left| f \left( x_s(\zeta\!-\!T_s\!+\!p) \right) \right|  
\exp \left[ \int_\zeta^\tau  {\mathrm{Tr}} \left[ D f \left( x_s(\xi\!-\!T_s\!+\!p) \right) \right] \d \xi \right] 
\d \zeta }{\left| f \left( x_s(\tau\!-\!T_s\!+\!p) \right) \right|^2} \d \tau \right| \nonumber \\
\label{eq:stableerrorboundparallel}
\end{eqnarray}
\end{small}
for  $ (p,t,\eps) \in \Xi_s \times (0,\eps_0) $, and (subject to being in $ \Xi_s $) obeys the 
limiting behaviour
\begin{small}
\begin{equation}
\lim_{t \rightarrow \infty} \left| e_s^\|(p,t) \right| = 0 \, , \, 
\lim_{p \rightarrow \infty} \left| e_s^\|(p,t) \right| \le \eps^2 \frac{\left( 2 C_g C_s + C_s^2 C_f \right)\left( \lambda_u + 2 C_f \right)}{- 2
 \lambda_s \lambda_u}  \left| 1 - e^{\lambda_s t} \right| \, .
\label{eq:stableerrorboundparallelinfinity} 
\end{equation}
\end{small}
\label{theorem:stable}
\end{theorem}

\proof
See Section~\ref{sec:stableproof}.
\proofend

Theorem~\ref{theorem:stable} provides a precise statement on why $ e_s $ is 
$ {\mathcal O}(\eps^2) $ for $ (p,t,\eps) \in \Xi_s\times
(0,\eps_0) $.  It should 
be noted that the improper integral in (\ref{eq:stableerrorboundperp}) is convergent, since 
as shown in the proof the integrand exhibits exponential decay.  Consequently, so
is the interior integral in (\ref{eq:stableerrorboundparallel}). The limiting
behaviour in (\ref{eq:stableerrorboundperpinfinity}) indicates how the perpendicular
component of the restricted
stable manifold error remains bounded in the limits as time goes to infinity, or in
each time-slice as the foot of the manifold (i.e., the hyperbolic trajectory $ a_\eps(t) $)
is approached.  The fact that the tangential component of the error approaches zero 
as $ t \rightarrow \infty $ is a consequence of the {\em restricted} nature of the stable
manifold.  As $ t \rightarrow \infty $, the length of the restricted stable manifold in each
time-section $ t $ goes to zero.  All points on these
one-dimensional curves---corresponding to all relevant $ p $ values---collapse together 
in the tangential direction, and as a consequence there is no error in this direction as $ t \rightarrow \infty $.
Put another way, both $ x_s^\eps $ and $ \tilde{x}_s^\eps $ undergo exponentially contracting behaviour
in the form $ e^{\lambda_s t} $ in the tangential direction, and so this is no surprise.

The detailed derivation of all this result is given in Section~\ref{sec:stableproof}, with a numerical example 
demonstrating the accuracy of the control strategy given in
Section~\ref{sec:taylorgreen}.

\section{Controlling unstable manifold}
\label{sec:unstablecontrol}

We now focus on determining the control velocity $ g $ in controlling the unstable manifold to have user-specified
behaviour.  The results are analogous to those of the stable manifold but require careful statement since there is
no requirement for the unstable manifold to have any relationship to the stable one.

Let $ T_u > 0 $ be a time-value before which the restricted
unstable manifold is to be quantified.  We represent 
the restricted two-dimensional
unstable manifold of (\ref{eq:pertaug}) when $ \eps = 0 $ by
\begin{equation}
\Gamma_u := \left\{ \left( x_u(t-T_u+p), t \right)  ~: ~ (p,t) \in (-\infty, U] \times 
(-\infty, T_u] \right\} \, ,
\label{eq:unpmanifoldu}
\end{equation}
in which $ x_u \left( \centerdot \right) $ is the trajectory lying along the unstable manifold.
The restricted unstable manifold which we desire to achieve in the $ \eps \ne 0 $ system will be represented by
\begin{equation}
\Gamma_u^\eps := \left\{ \left( x_u^\eps(p,t), t \right)  ~: ~ (p,t) \in (-\infty,U] \times (-\infty,T_u]  \right\} \, , 
\label{eq:pertmanifoldu} 
\end{equation}
for which we impose the conditions:

\begin{hypothesis}[Unstable manifold requirements]
For each $ t \le T_u $,  the quantity 
$ \left\{ x_u^\eps(p,t)~:~p \le U \right\} $ is a curve in 
$ \Omega $.  These restricted unstable manifold curves satisfy the following conditions.
\begin{itemize}
\item[(a)]  [Smoothness]  There exists a constant $ K_u > 0 $ such that for all
$ (p,t) \in (-\infty,U] \times (-\infty,T_u] $, and all $ \eps \in (0,\eps_0) $,
\begin{small}
\begin{equation}
\hspace*{-1cm}  \left|  x_u^\eps(p,t) \right| + \left| \frac{\partial}{\partial t} x_u^\eps(p,t) \right| +
\left|  \frac{\partial}{\partial p} x_u^\eps(p,t) \right| + \left| \frac{\partial}{\partial \eps} x_u^\eps(p,t) \right| + \left| \frac{\partial^2}{\partial \eps^2} x_u^\eps(p,t) \right|
+ \left| \frac{\partial^3}{\partial \eps^3} x_u^\eps(p,t) \right| < K_u \, .
\label{eq:xusmooth} 
\end{equation}
\end{small}
\item[(b)] [Closeness]  There exists a constant $ C_u > 0 $ such that for all
$ (p,t) \in (-\infty,U] \times (-\infty,T_u] $, 
\begin{equation}
\left| x_u^\eps (p,t) - x_u(t-T_u+p) \right| + \left| 
\frac{\partial}{\partial p} \left(x_u^\eps (p,t) - x_u(t-T_u+p) \right) \right|
\le C_u \, \eps \, . 
\label{eq:xuclose} 
\end{equation}
\item[(c)] [Limit]  For each $ t \le T_u $,  $ 
\displaystyle \lim_{p \rightarrow -\infty} x_u^\eps(p,t) $ is well defined.
\item[(d)] [Mappability] For each $ t \le T_u $, there exist intervals
$ \left( U_1(t), U_2(t) \right] $ and $ \left( U_1^\eps(t), U_2^\eps(t) \right]  $ ---both
of which are contained in $ (\infty, U]  $---  and 
 a scalar function $ r_u(\centerdot,t) $ defined on $ \left( U_1(t), U_2(t) \right] $ 
 which satisfies 
 \begin{equation}
 x_u^\eps(p^\eps,t) = x_u(t-T_u+p) + r_u(p,t)  \frac{J f \left( x_u(t-T_u+p) \right) }{ \left|  J f \left( x_u(t-T_u+p) \right)  \right|} \, ,
 \label{eq:mappableu}
 \end{equation}
 such that the mapping $ p \rightarrow p^\eps $ 
 from $ \left( U_1(t), U_2(t) \right] $ to $ \left( U_1^\eps(t), U_2^\eps(t) \right] $ 
 defined through (\ref{eq:mappableu})
 is a diffeomorphism.
 \item[(e)] [Congruence at time zero]   The $ p $ parameters in the time-slice $ t = 0 $ between
the unperturbed and the required restricted stable manifold curves match up, i.e., 
for all $ p \in \left( U_1(0), U_2(0) \right] $, 
\begin{equation}
\left[f \left( x_u(-T_u+p) \right) \right]^T  \left[ x_u^\eps(p,0) - x_u(-T_u+p) \right] = 0 \, .
\label{eq:xucongruence}
\end{equation}
\end{itemize}
\label{hyp:xu}
\end{hypothesis}

The set of $ (p,t) $ for which control is to be achieved is restricted to the set
\begin{equation}
\Xi_u^ := \left\{ \left(p, t \right)  ~: ~ t \le T_u~{\mathrm{and}}~  U_1(t) < p \le U_2(t)  \right\}  \,  ,
\label{eq:ptsetu} 
\end{equation}
where the largest intervals $ \left( U_1(t), U_2(t) \right] $ is chosen for each $ t $ in order to fulfil the 
mappability condition of Hypothesis~\ref{hyp:xu}.  
Now, for a prescribed restricted
unstable manifold $ x_u^\eps(p,t) $ we define the functions
\begin{equation}
M_u^\eps(p,t) := \left[ J f \left( x_u(t-T_u+p) \right) \right]^T \frac{x_u^\eps (p,t) - x_u(t-T_u+p)}{\eps}
\label{eq:mu}
\end{equation}
and
\begin{eqnarray}
B_u^\eps(p,t) :=  \left[ f \left( x_u(t-T_u+p) \right) \right]^T \frac{x_u^\eps (p,t) - x_u(t-T_u+p)}{\eps} \, ,
\label{eq:bu}
\end{eqnarray}
valid for $ (p,t,\eps) \in \Xi_u \times (0,\eps_0) $.

\begin{definition}[Control velocity for unstable manifold]
The control velocity $ g $ satisfies the smoothness conditions of Hypothesis~\ref{hyp:fgbound}, and moreover is
specified in normal and tangential components on the original unstable manifold by
\begin{eqnarray}
g^\perp \left( x_u(t-T_u+p),t \right)  & := & \frac{\left[ J f \left( x_u(t-T_u+p) \right) \right]^T}{
\left|  f \left( x_u(t-T_u+p) \right) \right|} g \left( x_u(t-T_u+p), t \right) \nonumber \\ 
& = & \frac{\frac{\partial M_u^\eps}{\partial t} (p,t) - {\mathrm{Tr}} \, 
\left( D f  \right) M_u^\eps(p,t)}{\left| f  \right|}
\label{eq:guperp}
\end{eqnarray}
and
\begin{eqnarray}
g^{\|} \left( x_u(t-T_u+p),t \right) & := & 
\frac{\left[ f \left( x_u(t-T_u+p) \right) \right]^T}{
\left|  f \left( x_u(t-T_u+p) \right) \right|} g \left( x_u(t-T_u+p), t \right) \nonumber \\
& & \hspace*{-2.5cm} = 
\frac{ \left| f \right|^2 \frac{\partial B_u^\eps(p,t)}{\partial t}
- f^T \left[ (D f) + (D f)^T \right] \left[ J f M_u^\eps(p,t) + f B_u^\eps(p,t)\right]}{ \left| f \right|^3} 
\label{eq:guparallel} 
\end{eqnarray}
in which $ f $ and $ D f $ in the above expressions are evaluated at $ x_u(t-T_u+p) $.
\label{def:unstable}
\end{definition}

Using the control velocity as defined in Def.~\ref{def:unstable} results in the required nonautonomous
unstable manifold to leading-order.
To characterise the resulting error, we define the true unstable manifold by
\begin{equation}
\tilde{\Gamma}_u^\eps := \left\{ \left( \tilde{x}_u^\eps(p,t), t \right)  ~: ~ (p,t) \in \Xi_u \right\} \, , 
\label{eq:pertmanifoldutilde} 
\end{equation}
rather than (\ref{eq:pertmanifoldu}), in which $ \tilde{x}_u^\eps(p,t) $ 
is an exact trajectory of (\ref{eq:pert}) which lies on the associated true perturbed manifold 
$ \tilde{\Gamma}_u^\eps $.  Analogous to the
congruence condition (\ref{eq:xucongruence}) at time zero for the desired unstable manifold, we require that 
\begin{equation}
\left[f \left( x_u(-T_u+p) \right) \right]^T  \left[ \tilde{x}_u^\eps(p,0) - x_u(-T_u+p) \right] = 0
\label{eq:xtildeucongruence}
\end{equation}
for the {\em true} unstable manifold.  The error in the restricted unstable manifold at
time $ t $ and parameter value $ p $ by $ e_u(p,t) $ is defined through
\begin{equation}
\tilde{x}_u^\eps(p,t) =  x_u^\eps(p,t) + e_u(p,t) \, .
\label{eq:xutilde}
\end{equation}

\begin{theorem}[Error in unstable manifold]
Assume the control velocity $ g $ satisfies Def.~\ref{def:unstable}, and define
\begin{equation}
e_u^\perp(p,t) := \frac{J f \left( x_u(t-T_u+p) \right)}{\left| f \left( x_u(t-T_u+p) \right) \right|}
e_u(p,t)
~~{\mathrm{and}}~~
e_u^\|(p,t) := \frac{ f \left( x_u(t-T_u+p) \right)}{\left| f \left( x_u(t-T_u+p) \right) \right|} e_u(p,t)
\, .
\label{eq:errorunstabledef}
\end{equation}
The normal component is bounded by
\begin{small}
\begin{equation}
 \left| e_u^\perp(p,t) \right| \le  \left[ C_u C_g + \frac{C_u^2 C_f}{2} \right] \eps^2 \frac{\int_{-\infty}^t 
\left| f \left( x_u(\tau\!-\!T_u\!+\!p) \right) \right|  
\exp \left[ \int_\tau^t  {\mathrm{Tr}} \left[ D f \left( x_u(\xi\!-\!T_u\!+\!p) \right) \right] \d \xi \right] 
\d \tau}{\left| f \left( x_u(t\!-\!T_u\!+\!p) \right) \right|} \, ,
\label{eq:unstableerrorboundperp}
\end{equation}
\end{small}
 for $ (p,t,\eps) \in \Xi_u \times (0,\eps_0) $, and  satisfies
\begin{equation}
\lim_{t \rightarrow -\infty} \left| e_u^\perp(p,t) \right| \le \frac{ \left[ C_u C_g + \frac{C_u^2 C_f}{2} \right] \eps^2}{\lambda_u} \, , \, 
\lim_{p \rightarrow -\infty} \left| e_u^\perp(p,t) \right| \le - \frac{ \left[ C_u C_g + \frac{C_u^2 C_f}{2} \right] \eps^2}{\lambda_s} \, ,
\label{eq:unstableerrorboundperpinfinity}
\end{equation}
as long as these limits can be taken within the domain $ \Xi_u $.
The tangential component of the error is
bounded by
\begin{small}
\begin{eqnarray}
\left| e_u^\|(p,t) \right| & \le &  \eps^2 \left( C_u C_g + \frac{ C_u^2 C_f}{2} \right) \left| f \left( x_u(t\!-\!T_u\!+\!p) \right) \right| \nonumber \\
& & \hspace*{-2cm} \times \left| \int_0^t
\frac{ \left| f \left( x_u(\tau\!-\!T_u\!+\!p) \right) \right| + 2 C_f \int_{-\infty}^\tau 
\left| f \left( x_u(\zeta\!-\!T_u\!+\!p) \right) \right|  
\exp \left[ \int_\zeta^\tau  {\mathrm{Tr}} \left[ D f \left( x_u(\xi\!-\!T_u\!+\!p) \right) \right] \d \xi \right] 
\d \zeta }{\left| f \left( x_u(\tau-T_u+p) \right) \right|^2} \d \tau \right| \nonumber \\
\label{eq:unstableerrorboundparallel}
\end{eqnarray}
\end{small}
for  $ (p,t,\eps) \in \Xi_u \times (0,\eps_0) $, 
and (subject to being in $ \Xi_u $) obeys the  limiting behaviour
\begin{small}
\begin{equation}
\lim_{t \rightarrow -\infty} \left| e_u^\|(p,t) \right| = 0 \, , \, 
\lim_{p \rightarrow -\infty} \left| e_u^\|(p,t) \right| \le \eps^2 \frac{\left( 2 C_g C_u + C_u^2 C_f \right)\left( - \lambda_s + 2 C_f \right)}{- 2
 \lambda_s \lambda_u}  \left| 1 - e^{\lambda_u t} \right| \, .
\label{eq:unstableerrorboundparallelinfinity} 
\end{equation}
\end{small}
\label{theorem:unstable}
\end{theorem}

\proof
See Section~\ref{sec:unstableproof}. 
\proofend

\section{Proof of Theorem~\ref{theorem:stable}}
\label{sec:stableproof}

We begin by introducing the notation
\begin{equation}
y(t) := x_s \left( t - T_s + p \right) \, ,
\label{eq:ys}
\end{equation}
which will be frequently needed in what follows.
We first argue that $ e_s(p,t) $ is bounded for $ (p,t,\eps) \in \Xi_s\times (0,\eps_0) $.  This is since
\begin{eqnarray*}
\left| e_s(p,t) \right| & = & \left| \tilde{x}_s^\eps(p,t) - x_s^\eps(p,t) \right| \\
& = & \left|  \tilde{x}_s^\eps(p,t) - a_\eps(t) + a_\eps(t) - a + a - y(t) + 
y(t) - x_s^\eps(p,t) \right| \\
& \le & \left|  \tilde{x}_s^\eps(p,t) - a_\eps(t) \right| + 
\left| a_\eps(t) - a  \right| + \left| a - y(t) \right| +
\left| y(t) - x_s^\eps(p,t) \right|
\end{eqnarray*}
We note that the first term goes to zero as $ t \rightarrow \infty $, 
since $ \tilde{x}_s^\eps(p,t) $ is an exact solution to the perturbed equation
(\ref{eq:pert}) which lies on the stable manifold of $ a_\eps(t) $.  The $ p $
selects a particular trajectory on this stable manifold, and thus this limit holds
for any $ p \ge S $.  Similarly, since 
$ y(t) = x_s(t-T_s+p) $ is on the stable manifold of $ a $, the third term also goes to
zero as $ t \rightarrow \infty $. Thus, these two terms are bounded.
The term $ \left| a_\eps(t) - a
\right| \le \eps C $ for some constant $ C $ for $ t \in [T_s,\infty) $ since the hyperbolic
trajectory remains $ {\mathcal O}(\eps) $-close to the unperturbed one \cite{coppel,yi}.  
Finally, the term $ \left| y(t) - x_s^\eps(p,t) \right| \le C_s \eps $ by 
Hypothesis~\ref{hyp:xs}.   Therefore, $ \left| e_s(p,t) \right| $ is bounded.

In contrast to $ M_s^\eps(p,t) $ in (\ref{eq:ms}), we define on $ \Xi_s $ 
an ``$M_s^\eps $ with error'' function
\begin{small}
\begin{eqnarray}
\hspace*{-0.5cm} \hat{M}_s^\eps(p,t) & :=  & \left[ J f \left( x_s(t-T_s+p) \right) \right]^T \frac{
\tilde{x}_s^\eps (p,t)  - x_s(t-T_s+p)}{\eps} \nonumber \\
& = & \left[ J f \left( y(t) \right) \right]^T \frac{
\tilde{x}_s^\eps (p,t)  - y(t)}{\eps} 
=  \left[ J f \left( y(t) \right) \right]^T \frac{
\left[ x_s^\eps (p,t) + e_s(p,t) \right] - y(t)}{\eps} \, .
\label{eq:mse}
\end{eqnarray}
\end{small}
The smoothness assumptions on $ f $ and $ g $ (Hypothesis~\ref{hyp:fgbound}) ensure
that the trajectory $ \tilde{x}_s^\eps(p,t) $ of (\ref{eq:pert}) 
is differentiable in $ t $ for any $ p $, and so
differentiating $ \hat{M}_s^\eps $ with respect to $ t $
leads to 
\begin{small}
\begin{eqnarray}
\hspace*{-2.5cm}
\eps \frac{\partial \hat{M}_s^\eps}{\partial t} (p,t) & = & \left[ J f \left( y(t) \right) \right]^T \left[ \frac{\partial \left[x_s^\eps (p,t)+ e_s(p,t) \right]}{\partial t}  - \frac{\partial y(t)}{\partial t} \right]  + \left[ J \, \frac{\partial f \left( y(t)\right)}{\partial t} \right]^T \left[ 
x_s^\eps (p,t) + e_s(p,t) - y(t) \right] \nonumber \\
& = & \left[ J f \left( y(t) \right) \right]^T \left[ f \left( x_s^\eps(p,t) + e_s(p,t)\right)
+ \eps g \left( x_s^\eps(p,t) + e_s(p,t), t \right) - f \left( y(t) \right)  \right] \nonumber \\
& & + \left[ J \, D f \left( y(t) \right) \frac{\partial y(t)}{\partial t}  \right]^T \left[ 
x_s^\eps (p,t) - y(t) + e_s(p,t) \right] \nonumber \\
& = & \eps \left[ J f \left( y(t) \right) \right]^T  g \left( x_s^\eps(p,t) + e_s(p,t), t \right) 
+ \left[ J f \left( y(t) \right) \right]^T \left[ f \left( x_s^\eps(p,t) + e_s(p,t)\right)
- f \left( y(t) \right) \right] \nonumber \\
& & + \left[ J \, D f \left( y(t) \right) f \left( y(t) \right)  \right]^T \left[ 
x_s^\eps (p,t) - y(t) \right] + \left[ J \, D f \left( y(t) \right) f \left( y(t) \right)  \right]^T e_s(p,t)
\label{eq:mstemp}
\end{eqnarray}
\end{small}
In the above calculations, the facts that $ x_s^\eps(p,t) + e_s(p,t) $ is an exact 
solution to the nonautonomous equation (\ref{eq:pert}), and $ y(t) = x_s(t-T_s+p) $
similarly satisfies the autonomous equation (\ref{eq:unp}) have been used.
 We note from Taylor's theorem that
\begin{eqnarray}
f \left( x_s^\eps(p,t) + e_s(p,t) \right) & =  & f \left( y(t) \right) + D f \left( y(t) \right)
\left( x_s^\eps(p,t) + e_s(p,t)- y(t) \right) \nonumber \\
& & \hspace*{-2.5cm} + \frac{1}{2} \left( x_s^\eps(p,t) + e_s(p,t) - y(t) \right)^T D^2 f \left( \xi_1 \right) 
\left( x_s^\eps(p,t) + e_s(p,t) - y(t) \right) \, .
 \label{eq:fexpand}
\end{eqnarray}
and that
\begin{equation}
g \left( x_s^\eps(p,t) + e_s(p,t) ,t \right)  =  g \left( y(t), t \right) + D g \left( \xi_2, t \right) 
\left( x_s^\eps(p,t) + e_s(p,t) - y(t) \right)
\label{eq:gexpand}
\end{equation}
for some points $ \xi_{1,2} \in \Omega $.
We substitute these expansions into (\ref{eq:mstemp}) and divide
by $ \eps $,
thereby arriving at
\begin{small}
\begin{eqnarray*}
\frac{\partial \hat{M}_s^\eps}{\partial t} (p,t) & = & 
\left[ J f \left( y(t) \right) \right]^T  g \left( y(t), t \right)  
+ \left[ J f \left( y(t) \right) \right]^T D f \left( y(t) \right) 
\frac{x_s^\eps(p,t) + e_s(p,t) - y(t)}{\eps} \\
& & + \left[ J \, D f \left( y(t) \right) f \left( y(t) \right)  \right]^T 
\frac{x_s^\eps (p,t) - y(t)}{\eps}
+ \left[ J \, D f \left( y(t) \right) f \left( y(t) \right)  \right]^T 
\frac{e_s(p,t)}{\eps} \\
& & +  \left[ J f \left( y(t) \right) \right]^T E_s(p,t) \, . 
\end{eqnarray*}
\end{small}
Here $ E_s(p,t) $ is a higher-order term satisfying 
\begin{equation}
\left| E_s(p,t) \right| \le \eps \left[ \tilde{C}_s C_g + \frac{\tilde{C}_s^2 C_f}{2} \right] 
\label{eq:Esbound}
\end{equation}
using (\ref{eq:closenesss}) and 
the bounds in Hypotheses~\ref{hyp:fgbound} and \ref{hyp:xs}, 
valid for $ (p,t,\eps) \in \Xi_s \times (0,\eps_0) $.
Using the easily verifiable identity $ \left[ J b \right]^T A + \left[ J A b \right]^T =
\left( {\mathrm{Tr}} \, A \right) \left[ J b \right]^T $ for $ 2 \times 1 $ vectors $ b $ and
$ 2 \times 2 $ matrices $ A $, we get $ 
\left[ J f \right]^T (D f) + \left[ J (D f) f 
\right]^T = {\mathrm{Tr}} \, \left( D f \right) \left[ J f \right]^T $, and hence
\begin{eqnarray*}
\frac{\partial \hat{M}_s^\eps}{\partial t} (p,t) & = & \left[ J f \left( y(t) \right) \right]^T  g \left( y(t), t \right) 
+  {\mathrm{Tr}} \left[ D f \left( y(t) \right) \right] \left[ J f \left( y(t) \right)
\right]^T \frac{ x_s^\eps(p,t) - y(t)}{\eps} \\
& & +  {\mathrm{Tr}} \left[ D f \left( y(t) \right) \right] \left[ J f \left( y(t) \right)
\right]^T
\frac{e_s(p,t)}{\eps} 
+ \left[ J f \left( y(t) \right) \right]^T  E_s(p,t) \, .
\end{eqnarray*}
Now, noting the definition of $ M_s^\eps(p,t) $ in comparison to $ \hat{M}_s^\eps(p,t) $, the above can be written as
\begin{small}
\begin{eqnarray}
& & \frac{\partial M_s^\eps}{\partial t}(p,t) + \frac{1}{\eps} \frac{\partial}{\partial t} \left\{
\left[ J f \left( y(t) \right) \right]^T e_s(p,t) \right\} = \left[ J f \left( y(t) \right) \right]^T  E_s(p,t) 
+ \left[ J f \left( y(t) \right) \right]^T  g \left( y(t), t \right) \nonumber \\
& & + {\mathrm{Tr}} \left[ D f \left( y(t) \right) \right] M_s^\eps(p,t) 
+ {\mathrm{Tr}} \left[ D f \left( y(t) \right) \right] \left[ J f \left( y(t) \right)
\right]^T
\frac{e_s(p,t)}{\eps}  \, .
\label{eq:mstemp2}
\end{eqnarray}
\end{small}
The intuition now is that we would like $ e_s $ to be $ {\mathcal O}(\eps^2) $, which
is yet to be established.   So we choose what we intend to be $ {\mathcal O}(\eps^0) $ terms above to be zero, that is, we
set
\[
\frac{\partial M_s^\eps}{\partial t}(p,t) =  \left[ J f \left( y(t) \right) \right]^T  
g \left( y(t), t \right) + {\mathrm{Tr}} \left[ D f \left( y(t) \right) \right] M_s^\eps(p,t) \, .
\]
Under this condition, we note that 
\[
g^\perp \left( y(t),t \right) :=  \frac{\left[ J f \left( y(t) \right) \right]^T}{
\left| f \left( y(t) \right) \right|} g \left( y(t), t \right) = \frac{\frac{\partial M_s^\eps}{\partial t} (p,t) - {\mathrm{Tr}} \left[ D f \left( y(t) \right) \right] M_s^\eps(p,t)}{\left| f \left( y(t) \right) \right|} \, , 
\]
which is exactly the control strategy defined in (\ref{eq:gsperp}).  Setting the normal
control velocity to equal this means that the remaining terms in (\ref{eq:mstemp2}) must
also be zero, that is
\[
 \frac{\partial}{\partial t} \left\{
\left[ J f \left( y(t) \right) \right]^T e_s(p,t) \right\} -   {\mathrm{Tr}} \left[ D f \left( y(t) \right) \right] \left[ J f \left( y(t) \right) \right]^T e_s(p,t) 
= \eps 
\left[ J f \left( y(t) \right) \right]^T  E_s(p,t) \, .
\]
Recalling the definition of $ y(t) $ in (\ref{eq:ys}), we multiply
through by the integrating factor 
\begin{equation}
\mu(p,t) := \exp \left[ \int_t^0  {\mathrm{Tr}} \left[ D f \left( y(\xi) \right) \right] \d \xi \right] \, , 
\label{eq:integratingfactor}
\end{equation}
giving the expression
\[
\frac{\partial}{\partial t} \left\{ \mu(p,t)
\left[ J f \left( y(t) \right) \right]^T e_s(p,t) \right\} = \eps \mu(p,t) 
\left[ J f \left( y(t) \right) \right]^T  E_s(p,t)
\]
which we integrate from a general $ t $ value to a large value $ L $ to obtain
\begin{small}
\begin{equation}
\mu(p,L) \left[ J f \left( y(L) \right) \right]^T e_s(p,L) 
- \mu(p,t) \left[ J f \left( y(t) \right) \right]^T e_s(p,t) 
=  
\eps \int_t^L \mu(p,\tau) 
\left[ J f \left( y(\tau) \right) \right]^T E_s(p,\tau)    
\d \tau \, .
\label{eq:stableerrortemp}
\end{equation}
\end{small}
We plan to take the limit $ L \rightarrow \infty $ in (\ref{eq:stableerrortemp}),
but first need to argue that this limit is defined.  Now
\begin{eqnarray*}
\left| \mu(p,L) \left[ J f \left( y(L) \right) \right]^T \right| 
& = & e^{\int_L^0  {\mathrm{Tr}} \left[ D f \left( y(\xi) \right) \right] \d \xi}
 \left| f \left( y(L) \right) \right| \\
& \rightarrow & e^{\int_L^0  \left( \lambda_s + \lambda_u \right) \d \xi}
  \, K \, e^{\lambda_s(L-T_s+p)} \\
& = & K e^{-(\lambda_u+\lambda_s)L} e^{\lambda_s(L-T_s+p)} \\
& = & K e^{\lambda_s(p-T_s)} e^{-\lambda_u L}
\end{eqnarray*}
where we have used the facts that $ {\mathrm{Tr}} \, \left( D f \right) $ approaches the sum of the
eigenvalues at $ a $ as its argument approaches $ a $, and that $ | f | $ has
exponential decay with rate $ \lambda_s $ as its argument approaches $ a $ along
the stable manifold.  Here, $ K $ is some constant, and since $ p \ge S $ and
$ \lambda_s < 0 $, the first exponential term is bounded by $ e^{\lambda_s(S-T_s)} $.
Thus, the quantity $ \left| \mu(p,L) \left[ J f \left( y(L) \right) \right]^T \right| $
decays exponentially in $ L $ with rate $ - \lambda_u $ as $ L \rightarrow \infty $.
We have at the beginning of this section argued that $ e_s(p,t) $ is bounded, and thus
when taking the limit $ L \rightarrow \infty $ in (\ref{eq:stableerrortemp}), the first
term on the left-hand side disappears.  On the other hand, the boundedness of 
$ E_s(p,t) $ given in (\ref{eq:Esbound}) in conjunction with the fact that the other
terms inside the integrand have $ e^{-\lambda_u \tau} $ behaviour (by the same
argument used above) implies that the 
improper integral on the right converges.  Thus we get
\begin{equation}
- \mu(p,t) \left[ J f \left( y(t) \right) \right]^T e_s(p,t) \nonumber \\
= 
\eps \int_t^\infty \mu(p,\tau) 
\left[ J f \left( y(\tau) \right) \right]^T E_s(p,\tau)    
\d \tau \, .
\label{eq:stableerrortemp2}
\end{equation}
Now we note from (\ref{eq:Esbound}) that
\begin{small}
\[
- \eps \left[ \tilde{C}_s C_g + \frac{ \tilde{C}_s^2 C_f}{2} \right] \left| f \left( y(\tau) \right) \right|  \le  \left[ J f \left( 
y(\tau) \right) \right]^T E_s(p,t)  \le 
\eps \left[ \tilde{C}_s C_g + \frac{ \tilde{C}_s^2 C_f}{2} \right] \left| f \left( y(\tau) \right) \right| \, .
\]
\end{small}
Dividing (\ref{eq:stableerrortemp2}) by $ \mu(p,t) \left| f \left( y(\tau) \right) \right| $ and utilising the above bounds, we get
\[
 \left| e_s^\perp(p,t) \right| \le  \left[ \tilde{C}_s C_g +
\frac{\tilde{C}_s^2 C_f}{2} \right] \eps^2 \frac{\int_t^\infty 
\left| f \left( y(\tau) \right) \right|  
e^{\int_\tau^t  {\mathrm{Tr}} \left[ D f \left( y(\xi) \right) \right] \d \xi}
\d \tau}{\left| f \left( y(t) \right) \right|} \, ,
\]
which is a genuine bound since the integrand of the improper integral exhibits exponential decay, and hence the integral is bounded.
Now, from (\ref{eq:xsclose}) and (\ref{eq:closenesss}) we see that it is possible to choose
$ \tilde{C}_s $ such that $ \left| \tilde{C}_s - C_s \right| \rightarrow 0 $ as $ \eps 
\rightarrow 0 $, and hence for sufficiently small $ \eps_0 $ is it possible to replace 
$ \tilde{C}_s $ above with $ C_s $, which leads directly to (\ref{eq:stableerrorboundperp}).
Moreover, the value of $ \left| e_s^\perp(p,t) \right| $ is
bounded as $ t \rightarrow \infty $, which is seen by a L'H\^opital's rule application to
the above:
\begin{eqnarray*}
\lim_{t \rightarrow \infty} \left| e_s^\perp(p,t) \right| & \le  & \left[ C_s C_g +
\frac{C_s^2 C_f}{2} \right] \eps^2 \lim_{t \rightarrow \infty} \frac{- \left| f \left( y(t) \right) \right|}{\frac{\partial}{\partial t} 
\left| f \left( y(t) \right) \right|} \\
& = &  - \left[ C_s C_g +
\frac{C_s^2 C_f}{2} \right] \eps^2 \lim_{t \rightarrow \infty} \frac{1}{\frac{\partial}{\partial t} \ln \left| f \left( y(t) \right) \right|} \, .
\end{eqnarray*}
But since  $ \left| f \left( y(t) \right) \right| \sim e^{\lambda_s (t-T_s+p)} $,  the
limit above is $ 1/ \lambda_s $, and we obtain the result in (\ref{eq:stableerrorboundperpinfinity}).  The limit $ p \rightarrow \infty $ at each fixed $ t $
is easiest computed with the formal 
replacements $ \left| f \left( y(t) 
\right) \right| \sim e^{\lambda_s(t-T_s+p)} $ and $ {\mathrm{Tr}} \, D f \left( y(\xi)
\right) \sim \lambda_u + \lambda_s $.  Thus,
\begin{eqnarray*}
\lim_{p \rightarrow \infty} \left| e_s^\perp(p,t) \right| & \le  & \left[ C_s C_g +
\frac{C_s^2 C_f}{2} \right] \eps^2 
\int_t^\infty \frac{e^{\lambda_s(\tau-T_s+p)}}{e^{\lambda_s(t-T_s+p)}} 
\exp \left[ \int_\tau^t \left( \lambda_s + \lambda_u \right) \, \d \xi \right] \d \tau \\
& = & \left[ C_s C_g +
\frac{C_s^2 C_f}{2} \right] \eps^2 
e^{\lambda_u t} \int_t^\infty e^{- \lambda_u \tau} \, \d \tau = \left[ C_s C_g +
\frac{C_s^2 C_f}{2} \right] \eps^2 \frac{1}{\lambda_u} \, .
\end{eqnarray*}
Hence 
$ \left| e^\perp(p,t) \right| $
exhibits the limiting behaviour in (\ref{eq:stableerrorboundperpinfinity}).
 
To evaluate the velocity requirement in the direction tangential to the manifold, we
proceed analogously and define
\begin{eqnarray}
\hspace*{-0.5cm} \hat{B}_s^\eps(p,t) & :=  & \left[ f \left( x_s(t-T_s+p) \right) \right]^T 
\frac{ \tilde{x}_s^\eps (p,t) - x_s(t-T_s+p)}{\eps}
\nonumber \\
& = & \left[ f \left( y(t) \right) \right]^T 
\frac{ \tilde{x}_s^\eps (p,t) - y(t)}{\eps}
=  \left[ f \left( y(t) \right) \right]^T \frac{
\left[ x_s^\eps (p,t) + e_s(p,t) \right] - y(t)}{\eps}
\label{eq:bse}
\end{eqnarray}
which differs from $ B_s^\eps $ in (\ref{eq:bs}) through the inclusion of the error
term $ e_s(p,t) $.
Taking the $ t $-derivative of $ \hat{B}_s^\eps $ leads to 
\begin{small}
\begin{eqnarray*}
\eps \frac{\partial \hat{B}_s^\eps}{\partial t} (p,t) & = & \left[  f \left( y(t) \right) \right]^T \left[ \frac{\partial}{\partial t}\left[  x_s^\eps (p,t) + e_s(p,t) \right] - \frac{\partial y(t)}{\partial t} \right] + \left[  \frac{\partial f \left( y(t) \right)}{\partial t} \right]^T \left[ 
x_s^\eps (p,t) + e_s(p,t) - y(t) \right] \\
& = & \left[ f \left( y(t) \right) \right]^T \left[ f \left( x_s^\eps(p,t) + e_s(p,t)\right)
+ \eps g \left( x_s^\eps(p,t) + e_s(p,t), t \right) - f \left( y(t) \right) \right] \\
& & + \left[ D f \left( y(t) \right) \frac{\partial y(t)}{\partial t}  \right]^T \left[ 
x_s^\eps (p,t) + e_s(p,t) - y(t) \right]\\
& = & \eps \left[ f \left( y(t) \right) \right]^T  g \left( x_s^\eps(p,t) + e_s(p,t), t \right) 
+ \left[f \left( y(t) \right) \right]^T \left[ f \left( x_s^\eps(p,t) + e_s(p,t)\right)
- f \left( y(t) \right) \right] \\
& & + \left[ D f \left(y(t) \right) f \left( y(t) \right)  \right]^T \left[ 
x_s^\eps (p,t) + e_s(p,t) - y(t) \right] \, .\\
\end{eqnarray*}
\end{small}
Applying the expansions (\ref{eq:fexpand}) and (\ref{eq:gexpand}) and dividing by
$ \eps $ gives
\begin{eqnarray*}
\frac{\partial \hat{B}_s^\eps}{\partial t} (p,t) & = &  \left[ f \left( y(t) \right) \right]^T  g \left( y(t), t \right) 
+ \left[ f \left( y(t) \right) \right]^T E_s(p,t) \\
& & + \left\{ \left[ f \left( y(t) \right) \right]^T D f \left( y(t)  \right) + \left[ D f \left( y(t)  \right) 
 f \left( y(t)  \right) \right]^T \right\}
\frac{x_s^\eps(p,t) + e_s(p,t) - y(t) }{\eps} 
 \end{eqnarray*} 
in which $ E_s(p,t) $ takes the same meaning as before, and satisfies the 
bound (\ref{eq:Esbound}).  Therefore,
\begin{eqnarray}
& & \frac{\partial B_s^\eps}{\partial t}(p,t) + \frac{1}{\eps} \frac{\partial}{\partial t} \left\{
\left[ f \left( y(t)  \right) \right]^T e_s(p,t) \right\} = \left[ f \left(y(t)  \right) \right]^T  E_s(p,t)
+ \left[ f \left( y(t)  \right) \right]^T  g \left( y(t), t \right) \nonumber \\
& & +  \left. f^T \left[ D f + \left( D f \right)^T \right] \right|_{y(t)} 
\frac{x_s^\eps(p,t) - y(t) }{\eps}
+   \left. f^T \left[ D f + \left( D f \right)^T \right] \right|_{y(t)}
\frac{e_s(p,t)}{\eps}  \, .
\label{eq:bstemp}
\end{eqnarray}
Now we write
\begin{equation}
\frac{x_s^\eps(p,t) - y(t) }{\eps} =  \frac{ J f \left( y(t)  \right) }{
\left| f \left( y(t)  \right) \right|^2} M_s^\eps(p,t) + 
  \frac{ f \left( y(t)  \right) }{
\left| f \left( y(t)  \right) \right|^2} B_s^\eps(p,t) \, , 
\label{eq:split}
\end{equation}
which is possible by (\ref{eq:ms}) and (\ref{eq:bs}) since $ M_s^\eps / \left| f \right| $
and $ B_s^\eps / \left| f \right| $ are the projections of the vector on the left-hand side
of (\ref{eq:split}) into the orthogonal directions given by $ J f / \left| f \right| $ and 
$ f / \left| f \right| $ respectively.  Substituting into (\ref{eq:bstemp}) yields
\begin{small}
\begin{eqnarray}
\frac{\partial B_s^\eps}{\partial t}(p,t) + \frac{1}{\eps} \frac{\partial}{\partial t} \left\{
\left[ f \left( y(t)  \right) \right]^T e_s(p,t) \right\} & = & \left[ f \left( y(t)  \right) \right]^T  E_s(p,t) + \left[ f \left(y(t)  \right) \right]^T  g \left( y(t) , t \right) \nonumber \\
& & \hspace*{-4cm} +  \left. \frac{f^T \left[ D f + \left( D f \right)^T \right] J f}{\left| f \right|^2} \right|_{y(t) } 
M_s^\eps(p,t) +  \left. \frac{f^T \left[ D f + \left( D f \right)^T \right] f}{\left| f \right|^2} \right|_{y(t) } 
B_s^\eps(p,t) \nonumber \\
& & \hspace*{-4cm} +   \left. f^T \left[ D f + \left( D f \right)^T \right] \right|_{y(t) }
\frac{e_s(p,t)}{\eps}  \, .
\label{eq:bstemp2}
\end{eqnarray}
\end{small}
We select the terms we plan to be $ {\mathcal O}(\eps^0) $ above to be zero, giving
\begin{eqnarray*}
\hspace*{-1cm}
\left[ f \left( y(t)  \right) \right]^T  g \left( y(t) , t \right) & = &
\frac{\partial B_s^\eps}{\partial t} (p,t) -  \left. \frac{f^T \left[ D f + \left( D f \right)^T \right] J f}{\left| f \right|^2} \right|_{y(t) } 
M_s^\eps(p,t) \\
& & -  \left. \frac{f^T \left[ D f + \left( D f \right)^T \right] f}{\left| f \right|^2} \right|_{y(t) } 
B_s^\eps(p,t) 
\end{eqnarray*}
Thus
\begin{small}
\begin{eqnarray*}
g_s^{\|} \left( y(t) ,t \right) & := & \frac{\left[ f \left( y(t)  \right) \right]^T}{
\left|  f \left( y(t)  \right) \right|} g \left( y(t) , t\right) \nonumber \\
& = & \frac{1}{\left| f \left( y(t)  \right) \right|} 
\frac{\partial B_s^\eps}{\partial t} (p,t) -  \left. \frac{f^T \left[ D f + \left( D f \right)^T \right]}{\left| f \right|^3} \right|_{y(t) }
\left( J f M_s^\eps(p,t) + f B_s^\eps(p,t) \right)  \, , 
\end{eqnarray*}
\end{small}
which is the
tangential component of the control velocity required, as given in equation
(\ref{eq:gsparallel}).  Under this choice,  the remaining 
 terms in (\ref{eq:bstemp2}) must equal zero, and hence
 \begin{small}
\begin{eqnarray}
\hspace*{-1.2cm} \frac{\partial}{\partial t} \left\{
\left[ f \left( y(t)  \right) \right]^T e_s(p,t) \right\} & = & \eps \left[ f \left( y(t)  \right) \right]^T  E_s(p,t) +  \left. f^T \left[ D f + \left( D f \right)^T \right] \right|_{y(t) } e_s(p,t)
\, .
\label{eq:tangentialerrortemp}
\end{eqnarray}
\end{small}
We now write $ e_s(p,t) $ in terms of the orthogonal unit vectors $ J f / |f | $ and $
f /|f| $ as
\begin{small}
\[
e_s(p,t) = \left\{ \left[ J f \right]^T\left(y(t)  \right) e_s(p,t) \right\} \frac{J f \left( y(t) \right)}{\left| f \left( y(t)  \right) \right|^2}
 + \left\{  f^T\left( y(t)  \right) e_s(p,t) \right\} \frac{f \left( y(t)  \right)}{\left| f \left( y(t) \right) \right|^2}
 \, , 
\]
\end{small}
enabling (\ref{eq:tangentialerrortemp}) to be written as
\begin{small}
\[
\frac{\partial}{\partial t} \left[ f^T e_s(p,t) \right] = \eps f^T E_s(p,t) + 
\frac{f^T \left[ D f + \left( D f \right)^T \right] J f}{\left| f \right|^2} \left( J f \right)^T e_s + 
\frac{f^T \left[ D f + \left( D f \right)^T \right] f}{\left| f \right|^2} 
f^T e_s 
\]
\end{small}
where the argument $ y(t)  $ in each of the $ f $ terms has been suppressed
for convenience.  Thus we have the equation
\begin{small}
\begin{equation}
\frac{\partial}{\partial t} \left[ f^T e_s(p,t) \right] - \frac{f^T \left[ D f + \left( D f \right)^T \right] f}{\left| f \right|^2} \left[ f^T e_s(p,t) \right]
 = \eps f^T E_s(p,t) + \frac{f^T \left[ D f + \left( D f \right)^T \right] J f}{\left| f \right|} 
  \frac{\left( J f 
\right)^T}{\left| f \right|}  e_s 
\label{eq:tangentialerrorode}
\end{equation}
\end{small}
We will consider this a linear equation for $ f^T e_s $, since we will show that the 
right-hand side can be bounded.  The left-hand side can be simplified with the observation
\begin{eqnarray}
\frac{\partial}{\partial t} \left[ f \left( y(t) \right) \right] & = & D f \left( y(t) \right)
\frac{\partial}{\partial t} \left[ y(t) \right] \nonumber \\
& = & D f \left( y(t) \right) f \left( y(t) \right) \, , 
\label{eq:variational}
\end{eqnarray}
and that its transpose yields
\begin{equation}
\frac{\partial}{\partial t} \left[ f^T \left( y(t)  \right) \right] =
f^T \left( y(t) \right) \left[ D f \right]^T \left( y(t) \right) \, .
\label{eq:adjoint}
\end{equation}
Therefore, we note that
\[
\frac{\partial}{\partial t} \left[ f^T \left( y(t)  \right) f \left( y(t)  \right) 
\right] = f^T \left[ D f f \right] + \left[ f^T (D f)^T \right] f = f^T \left[ D f + (D f)^T \right] f
\]
Hence, 
\[
\frac{f^T \left[ D f + \left( D f \right)^T \right] f}{\left| f \right|^2} = \frac{1}{\left| f \right|^2}
\frac{\partial}{\partial t} \left[ f^T f \right] = \frac{1}{\left| f \right|^2} \frac{\partial}{\partial t}
\left| f \right|^2 = \frac{\partial}{\partial t} \ln \left| f \right|^2 \, , 
\]
and therefore (\ref{eq:tangentialerrorode}) can be written as
\begin{small}
\begin{equation}
\frac{\partial}{\partial t} \left[ f^T e_s(p,t) \right] - \frac{\partial}{\partial t} \left[ \ln \left| f \right|^2 \right] \left[ f^T e_s(p,t) \right]
 =  \eps f^T E_s(p,t) + \frac{f^T \left[ D f + \left( D f \right)^T \right] J f}{\left| f \right|} 
 e_s^\perp 
 \, .
\label{eq:tangentialerrorode2}
\end{equation}
\end{small}
where we have used the fact that $ e_s^\perp = \left( J f \right)^T e_s / \left| f \right| $.
Multiplying (\ref{eq:tangentialerrorode2}) through by the integrating factor
$ \left| f \left( y(t) \right) \right|^{-2} $, and integrating from $ 0 $ to a general $ t $
value yields
\begin{equation}
 \frac{f^T \left( y(t) \right)}{\left| f \left( y(t) \right) \right|^2}
 e_s(p,t) = 
 \int_0^t 
\frac{\eps f^T E_s(p,\tau) + \frac{f^T \left[ D f + \left( D f \right)^T \right] J f}{\left| f \right|} \Big|_{y(\tau)}
 e_s^\perp }{\left| f \left( y(\tau) \right) \right|^2} \d \tau \, ,
\label{eq:tangentialerror1}
\end{equation}
in which the congruence condition (\ref{eq:xscongruence}) has been used to get rid
of the boundary term at $ t = 0 $.
Now, we bound the integrand in (\ref{eq:tangentialerror1}) using 
(\ref{eq:stableerrorboundperp}), (\ref{eq:Esbound}) and 
Hypothesis~\ref{hyp:fgbound}, and with the understanding that $ \tilde{C}_s $ can be
replaced with $ C_s $ for suitably small $ \eps_0 $: 
\begin{eqnarray*}
& & \left| \frac{\eps f^T E_s(p,\tau) + \frac{f^T \left[ D f + \left( D f \right)^T \right] J f}{\left| f \right|} 
 e_s^\perp }{\left| f \right|^2} \right| \\
& \le & \frac{\left| f^T \right|}{\left| f \right|^2}  \left| \eps \left| E_s \right| + \left| D f + \left( D f \right)^T \right| \left|
 e_s^\perp \right| \right| \\
 & \le & \left| f \right|^{-1} \left| \eps^2 \left( C_s C_g + \frac{ C_s^2 C_f}{2} \right) + 2 C_f \left| e_s^\perp \right| \right| \\
 & \le & \left| f \right|^{-1} \eps^2 \left( C_s C_g + \frac{ C_s^2 C_f}{2} \right) \left[
 1 + 2 C_f \frac{\int_\tau^\infty 
\left| f \left( y(\zeta) \right) \right|  
e^{ \int_\zeta^\tau  {\mathrm{Tr}} \left[ D f \left( y(\xi) \right) \right] \d \xi} 
\d \zeta}{\left| f \left( y(\tau) \right) \right|} \right] \\
& =: & \left| f \right|^{-1} \eps^2 \left( C_s C_g + \frac{ C_s^2 C_f}{2} \right) H(p,\tau)
 \end{eqnarray*}
 which defines $ H $ as the term in the square brackets, and we note that $ H $ is
bounded in $ \tau $ since the $ \tau $-dependent quotient in $ H $ has a finite limit as established
in (\ref{eq:stableerrorboundperpinfinity}).
Therefore from (\ref{eq:tangentialerror1}),
\begin{small}
\begin{eqnarray*}
 \hspace*{-2cm} \left| e_s^\|(p,t) \right| & = & \left| \frac{f^T \left( y(t) \right)}{\left| f \left( y(t) \right) \right|}  e_s(p,t) \right| \nonumber \\
& \le &  \eps^2 \left( C_g C_s + \frac{ C_s^2 C_f}{2} \right) \left| f \left( y(t) \right) \right| \left| \int_0^t
\frac{H(p,\tau)}{\left| f \left( y(\tau) \right) \right|} \d \tau \right| \nonumber \\
& = & \eps^2 \left( C_g C_s + \frac{ C_s^2 C_f}{2} \right) \left| f \left( y(t) \right) \right|  \left| \int_0^t
\frac{ \left| f \left( y(\tau) \right) \right| + 2 C_f \int_\tau^\infty 
\left| f \left( y(\zeta) \right) \right|  
e^{\int_\zeta^\tau  {\mathrm{Tr}} \left[ D f \left( y(\xi) \right) \right] \d \xi}
\d \zeta }{\left| f \left( y(\tau) \right) \right|^2} \d \tau \right| \
\end{eqnarray*}
\end{small}
Writing in the $ \infty / \infty $ form, L'H\^opital's Rule can be used to show that the above goes to zero as $ t \rightarrow
\infty $:
\[
\lim_{t \rightarrow \infty} \frac{ \int_0^t 
\frac{H(p,\tau)}{\left| f \left( y(\tau) \right) \right|} \d \tau
}{ \left| f \left( y(t) \right) \right|^{-1}} =  \lim_{t \rightarrow \infty} \frac{  \frac{H(p,t)}{\left| f \left( y(t) \right) \right|} }{ - 
\left| f \left( y(t) \right) \right|^{-2}} = - \lim_{t \rightarrow \infty} H(p,t) \left| f \left( y(t) \right) \right| = 0 \, , 
\]
since $ H $ is bounded and $ \left| f \left( y(t) \right) \right| \rightarrow 
\left| f (a) \right| = 0 $.  
To take the $ p \rightarrow \infty $ limit, we proceed as
before and replace each term with its appropriate
limiting behaviour, and thus
\begin{small}
\begin{eqnarray*}
\hspace*{-2cm}  \lim_{p \rightarrow \infty} \left| e_s^\|(p,t) \right| & \le &  \eps^2 \left( C_g C_s + \frac{ C_s^2 C_f}{2} \right) e^{\lambda_s(t-T_s+p)} \left| \int_0^t \frac{e^{\lambda_s(\tau-T_s+p)} + 
 2 C_f \int_\tau^\infty e^{\lambda_s(\zeta-T_s+p)} e^{(\lambda_s+\lambda_u)(\tau-\zeta)}
 \d \zeta}{e^{2 \lambda_s(\tau-T_s+p)}} \, \d \tau \right|  \\
 & = &  \eps^2 \left( C_g C_s + \frac{ C_s^2 C_f}{2} \right) e^{\lambda_s(t-T_s+p)}
 \left| \int_0^t \frac{e^{\lambda_s(\tau-T_s+p)} + 
 2 C_f e^{\lambda_s(\tau-T_s+p)} e^{\lambda_u \tau} \int_\tau^\infty  e^{-\lambda_u \zeta}
 \d \zeta}{e^{2 \lambda_s(\tau-T_s+p)}} \, \d \tau \right| \\
 & = &  \eps^2 \left( C_g C_s + \frac{ C_s^2 C_f}{2} \right) e^{\lambda_s(t-T_s+p)}
 \left| \int_0^t e^{-\lambda_s(\tau-T_s+p)} \left[ 1 + 
 2 C_f e^{\lambda_u \tau} \int_\tau^\infty  e^{-\lambda_u \zeta}
 \d \zeta \right] \, \d \tau \right| \\
 & = &  \eps^2 \left( C_g C_s + \frac{ C_s^2 C_f}{2} \right) e^{\lambda_s t}
\left|  \int_0^t e^{-\lambda_s\tau} \left[ 1 + 
 \frac{2 C_f}{\lambda_u}  \right] \, \d \tau \right| \\
 & = & \eps^2 \left( C_g C_s + \frac{ C_s^2 C_f}{2} \right) \left( 1 + 
 \frac{2 C_f}{\lambda_u}  \right) e^{\lambda_s t}
 \left| \frac{ 1 - e^{-\lambda_s t}}{\lambda_s} \right| \\
 & = & \eps^2 \frac{\left( 2 C_g C_s + C_s^2 C_f \right)\left( \lambda_u + 2 C_f \right)}{2
 \left| \lambda_s \right| \lambda_u}  \left| e^{\lambda_s t} - 1 \right|
\end{eqnarray*}
\end{small}
which is  (\ref{eq:stableerrorboundparallelinfinity}) as required.
Thus, $ e^\| $ remains $ {\mathcal O}(\eps^2) $ just
as $ e^\perp $ does, implying that $ e_s(p,t) $ is $ {\mathcal O}(\eps^2) $ as desired.

\section{Proof of Theorem~\ref{theorem:unstable}}
\label{sec:unstableproof}

The proof is analogous to the stable manifold results, and requires the definitions
\begin{equation}
\hat{M}_u^\eps(p,t) :=  \left[ J f \left( x_u(t-T_u+p) \right) \right]^T \frac{
\left[ x_u^\eps (p,t) + e_u(p,t) \right] - x_u(t-T_u+p)}{\eps}
\label{eq:mue}
\end{equation}
and
\begin{equation}
\hat{B}_u^\eps(p,t) :=  \left[ f \left( x_u(t-T_u+p) \right) \right]^T \frac{
\left[ x_u^\eps (p,t) + e_u(p,t) \right] - x_u(t-T_u+p)}{\eps} \, . 
\label{eq:bue}
\end{equation}
The proof then proceeds exactly as in
Theorem~\ref{theorem:stable}, with the only substantive changes being that the subscript $ s $
(for stable) needs to be replaced with the subscript $ u $ (for unstable), and that  integration occurs from
$ - \infty $ to a general time as opposed to from a general time to $ + \infty $ when
working with the normal component of $ g $. Details
will not be provided.

\section{Taylor-Green flow example}
\label{sec:taylorgreen}

\begin{figure}[t]
\includegraphics[scale=0.8]{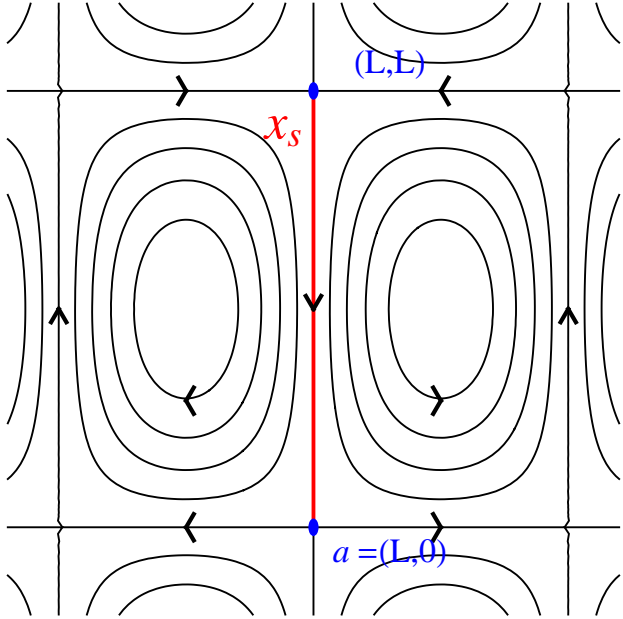}
\caption{The stable manifold branch of $ (L,0) $ in the Taylor-Green flow 
(\ref{eq:taylorgreen}) which is to be controlled.}
\label{fig:taylorgreen}
\end{figure}

We will present a short example to demonstrate the efficacy of the theoretical 
method, postponing an extensive numerical analysis to a future article.
Consider the Taylor-Green flow
\begin{equation}
\renewcommand{\arraystretch}{1.8}
\left. \begin{array}{c}
\dot{x} = - \pi U \sin \left( \frac{\pi x}{L} \right) \cos \left( \frac{\pi y}{L} \right) \\
\dot{y} = \pi U \cos \left( \frac{\pi x}{L} \right) \sin \left( \frac{\pi y}{L} \right)
\end{array}
\right\} \, ,
\label{eq:taylorgreen}
\end{equation}
in which $ U $ and $ L $ are positive parameters with dimensions of velocity and
length respectively.  This flow is equivalent to the steady limit of the popular 
double-gyre model \cite{shadden}.
The autonomous system (\ref{eq:taylorgreen}) possesses a heteroclinic trajectory from the fixed point $ (L,L) $ to that at
$ (L,0) $, which is given by
\[
\renewcommand{\arraystretch}{1.8}
x_{s,u}(t) = \left( \begin{array}{c} L \\ \frac{2 L}{\pi} \tan^{-1} e^{- \pi^2 U t/L} 
\end{array} \right) \, .
\]
The above notation has been used since this is the stable manifold of $ (L,0) $, but
is the unstable manifold of $ (L,L) $.  See Fig.~\ref{fig:taylorgreen}.  
Here, we will focus only on 
controlling the stable manifold $ x_s $ of the fixed point $ a \equiv (L,0) $.
Note in particular that since the manifold is downwards along the line $ x = L $, the perpendicular
and parallel components required in Def.~\ref{def:stable} relate exactly
to the $ x $ and $ -y $ directions at every point on the heteroclinic.  
As an example, we shall try to move this 
stable manifold to the nonautonomous location
\begin{equation}
\renewcommand{\arraystretch}{1.8}
x_s^\eps(p,t) = \left( \begin{array}{l} L \\ \frac{2L }{\pi} \tan^{-1} e^{- \pi^2 U (t-T_s+p)/L} 
\end{array} \right) + \eps L \left( \begin{array}{l} e^{-U p/L} \cos \frac{U (t-p)}{L} \\ 0
\end{array} \right) \, ,
\label{eq:taylorgreenstable}
\end{equation}
by introducing a control velocity $ g(x,y,t) $, which we shall in this case insist on
being incompressible to be consistent with the incompressibility of the Taylor-Green
flow.  The $ \eps = 0 $ version of (\ref{eq:taylorgreenstable}) is exactly 
$ x_s(t-T_s+p) $; we have built in the $ {\mathcal O}(\eps) $-closeness of the
desired manifold to the unperturbed one directly.
To determine the form of this curve in each time-slice $ t $, we can think (\ref{eq:taylorgreenstable}) at each fixed $ t $ value subject to $ t \ge T_s $.  This would
then be a parametric representation in terms of the parameter $ p \ge S $; we 
can take $ S_1(t) = S $ for all $ t $ and $ S_2(t) = \infty $ for this chosen form.  Thus,
the theory will work on $ \Xi_s = \left\{ (p,t)~:~ p \ge S ~{\mathrm{and}}~ t\ge T_s
\right\} $. The
beginning of this manifold in the time-slice $ t $, that is, the location of the hyperbolic
trajectory associated with the unperturbed saddle point $ (L,0) $, can be obtained
by taking the limit as $ p \rightarrow \infty $, which yields $ (L,0) $ for all $ t $.
We can indeed find the required stable manifold curve in each time-slice $ t $ by
eliminating $ p $ from the parametric equation (\ref{eq:taylorgreenstable}); since
\begin{equation}
y(p,t) = \frac{2L }{\pi} \tan^{-1} e^{- \pi^2 U (t-T_s+p)/L}
\label{eq:ypt}
\end{equation}
we have the relationship
\begin{equation}
p(y,t) = - \frac{L}{\pi^2 U} \ln \left( \tan \frac{\pi y}{2 L} \right) + T_s - t \, ,
\label{eq:pyt}
\end{equation}
and thus the stable manifold curve in each time-slice $ t $ in $ (x,y) $-coordinates
is
\begin{equation}
x = L \left\{ 1 + \eps \exp \left[ - \frac{U p(y,t)}{L} \right] \cos \frac{U \left( t - p(y,t) \right)}{L} 
\right\} \, ,
\label{eq:taylorgreenmanifold}
\end{equation}
subject to the restrictions $ t \ge T_s $ and $ p \ge S $.  The condition on $ p $ 
can be translated to
\begin{equation}
0 < y < y_m(t) := \frac{2 L}{\pi} \tan^{-1} \exp \left[ -\frac{\pi^2 U (t - T_s+S)}{L} \right] \, 
\label{eq:taylorgreenymax}
\end{equation}
where $ y_m(t) $ is the maximum value of $ y $ attainable in the time-slice $ t $. 
We observe that (\ref{eq:taylorgreenstable}) also satisfies  the 
congruence condition (\ref{eq:xscongruence}) since the $ {\mathcal O}(\eps) $
term in (\ref{eq:taylorgreenstable}) is in the $ x $-direction at $ t = 0 $, and is thus perpendicular to the
unperturbed stable manifold.  Now, in this case the components of the control 
$ g(x,y,t) $ we need are $ g^\perp $ (in the $ +x $ direction)
and $ g^\| $ (in the $ -y $ direction).  By utilising the requirements in Def.~\ref{def:stable} and doing the relevant algebra (not shown), we find that the control $ g $ 
needs to satisfy
\begin{small}
\[
g \left(L, \frac{2 L}{\pi} \tan^{-1} e^{-\pi^2 U(t-T_s+p)/L}, t \right) = -U e^{-Up / L} \left(
\begin{array}{l}
\sin  \frac{ U (t-p)}{L} + \pi^2 \tanh \frac{\pi^2 U (t - T_s + p)}{L}  \cos  \frac{ U (t-p)}{L} \\ 
0
\end{array} \right) \, .
\]
\end{small}
Any control velocity $ g(x,y,t) $ satisfying (\ref{eq:taylorgreencondition}) is appropriate.
We note that there are infinitely many ways to do this, since it is only the value of $ g $
on the stable manifold which needs to be specified.  We choose the following strategy to find one such $ g $.  By replacing $ p $ with (\ref{eq:pyt}), 
we realise that we have the relationship
\[
\tanh \left[ \frac{\pi^2 U (t-T_s+p)}{L} \right] = \cos \frac{\pi y}{L} 
\]
resulting in
\begin{equation}
\renewcommand{\arraystretch}{1.8}
g \left(L, y, t \right) = U \left(
\begin{array}{l}
- e^{-U p(y,t)/L} \left[  
\sin  \frac{ U (t-p(y,t))}{L} + \pi^2 \cos \frac{\pi y}{L}  \cos  \frac{ U (t-p(y,t))}{L} \right] \\ 
0
\end{array} \right) \, .
\label{eq:taylorgreencondition}
\end{equation}
Now, any form
for $ g(x,y,t) $ which is consistent with (\ref{eq:taylorgreencondition}) will result 
in our desired restricted stable manifold, correct to $ {\mathcal O}(\eps) $.  The
easiest option would be to
extend uniformly in $ x $, which can be seen to preserve incompressibility.
We will choose an alternative $ g $, determined by adding a divergence-free term to the above which 
yields zero when evaluated on $ x = L $, that is, we choose the control
\begin{equation}
\renewcommand{\arraystretch}{1.8}
g \left(x, y, t \right) = U \left(
\begin{array}{l}
-e^{-U p(y,t)/L} \left[  
\sin  \frac{ U (t-p(y,t))}{L} + \pi^2 \cos \frac{\pi y}{L}  \cos  \frac{ U (t-p(y,t))}{L} 
\right]  \\ 
 \sin \frac{\pi x}{L} \sin \frac{U t}{L}
\end{array} \right) \, .
\label{eq:taylorgreencontrol}
\end{equation}

\begin{figure}[t]
\begin{tabular}{c}
\includegraphics[width=10cm]{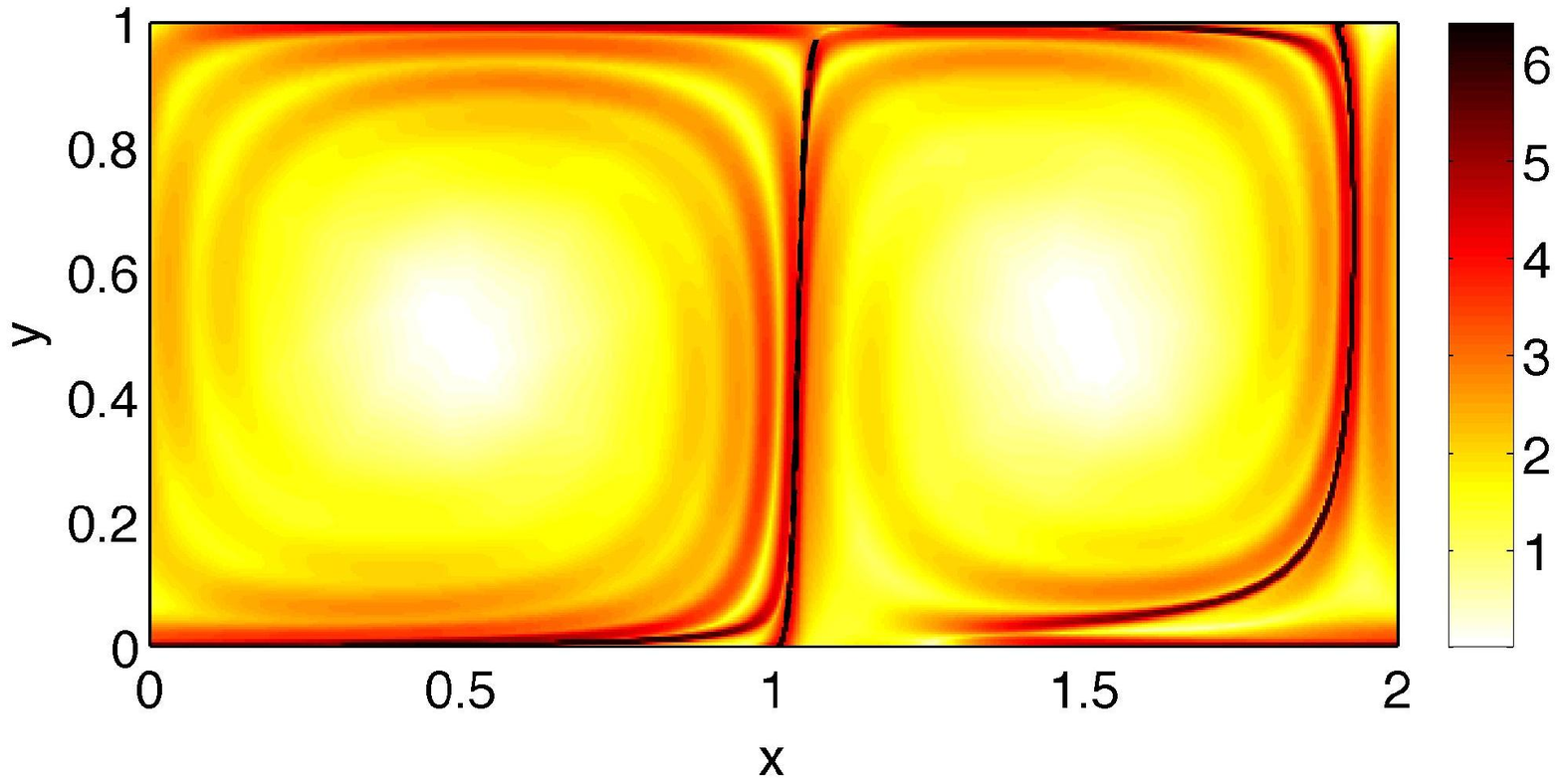}\\
\includegraphics[width=10cm]{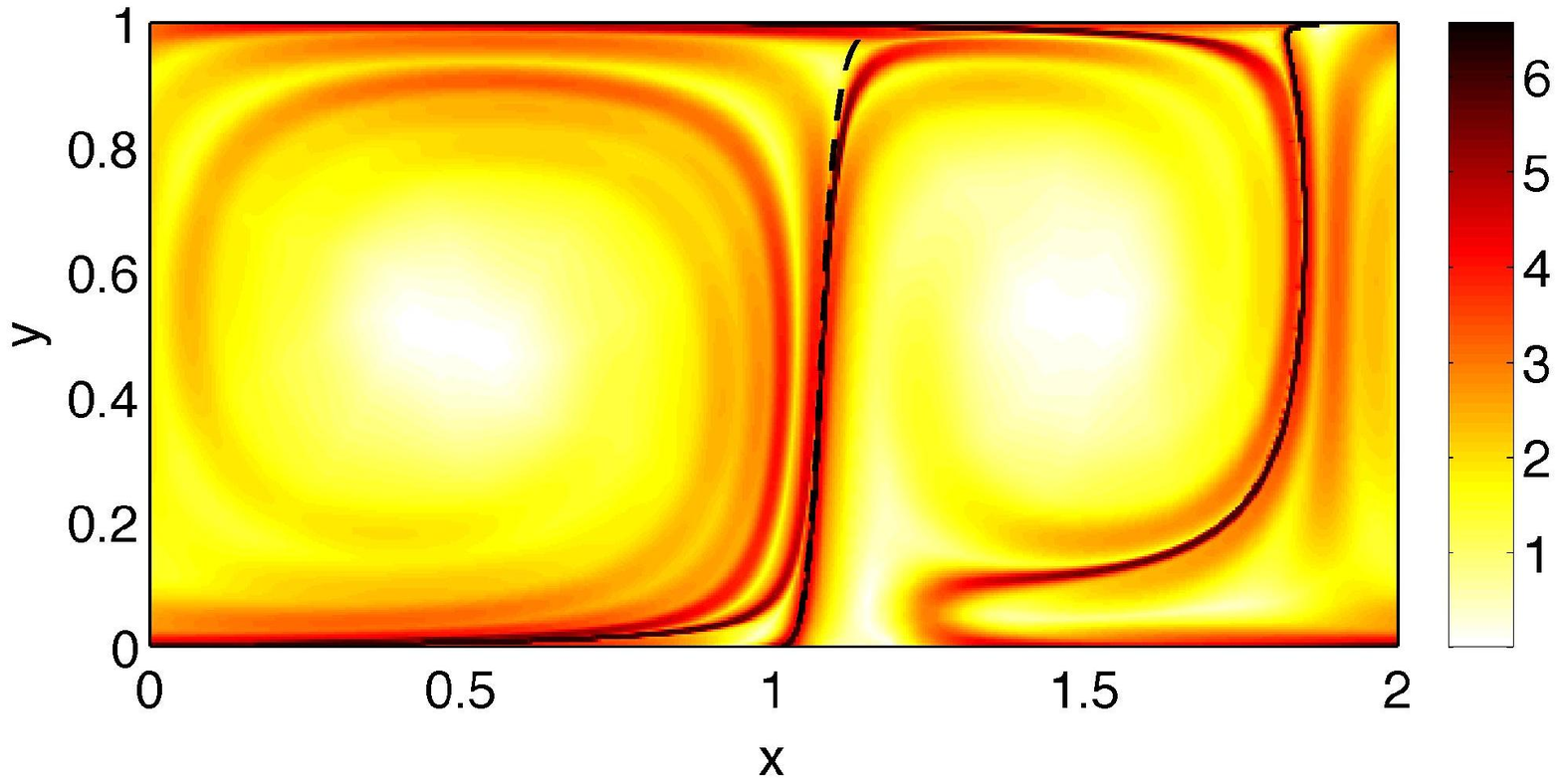}\\
\includegraphics[width=10cm]{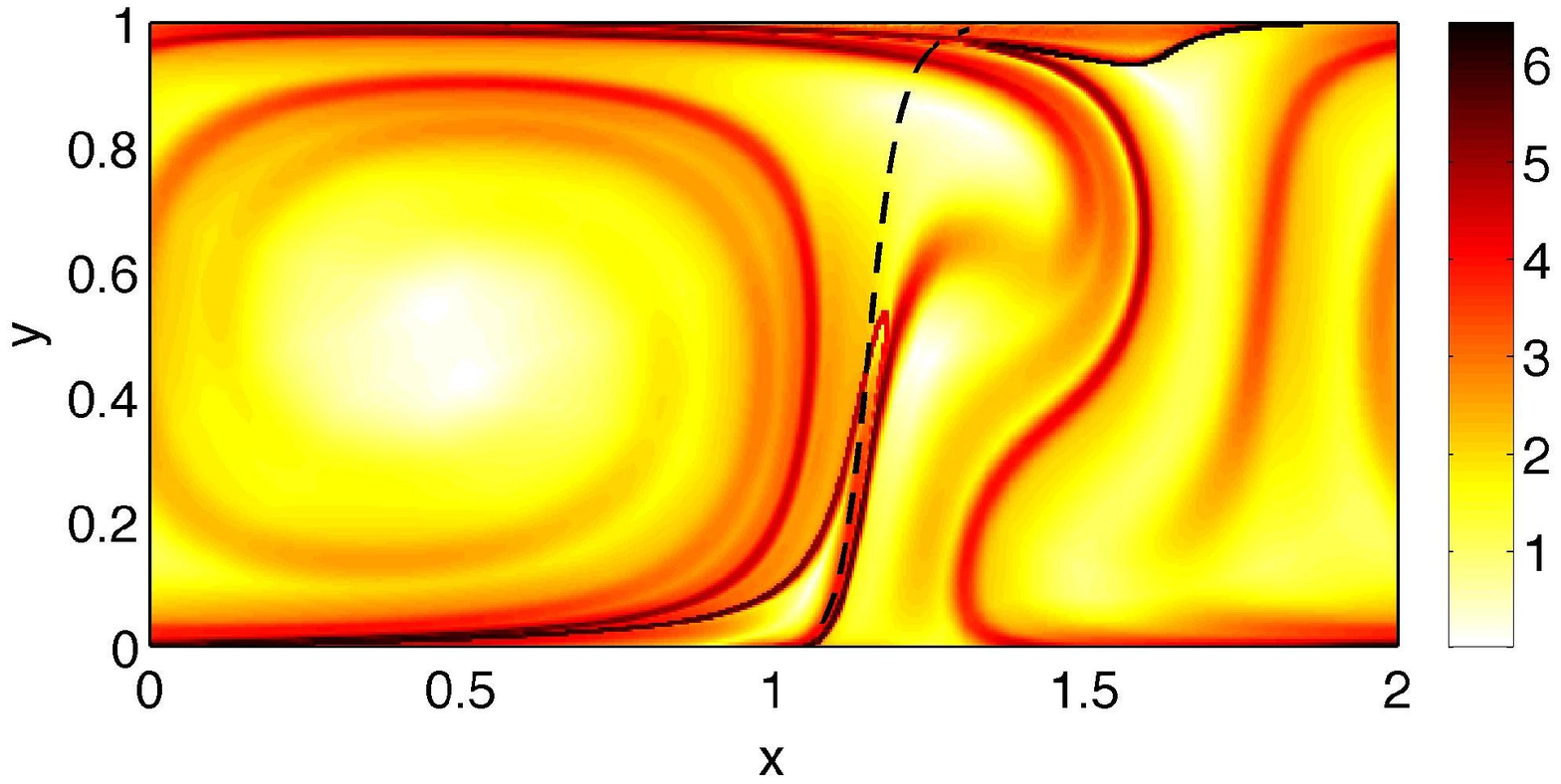}\\
\end{tabular}
\caption{Finite-time Lyapunov exponent fields at $t=-0.9$ for (\ref{eq:taylorgreencontrolled}) with $U=1$, $L=1$ and $T_s=-1$. The desired stable manifold (\ref{eq:taylorgreenmanifold}) is shown by the black dashed curve, and the panels are respectively for the choices $\eps=0.05, \, 0.1, \, 0.2$.}\label{fig:ftle_eps}
\end{figure}

Thus, the claim is that (\ref{eq:taylorgreenmanifold}) is the restricted stable manifold
of the system
\begin{small}
\begin{equation}
\renewcommand{\arraystretch}{1.8}
\left. \begin{array}{l}
\dot{x} = - \pi U \sin \left( \frac{\pi x}{L} \right) \cos \left( \frac{\pi y}{L} \right) - \eps U 
e^{-U p(y,t)/L} \left[  
\sin  \frac{ U (t-p(y,t))}{L} + \pi^2 \cos \frac{\pi y}{L}  \cos  \frac{ U (t-p(y,t))}{L} 
\right]  \\
\dot{y} = \pi U \cos \left( \frac{\pi x}{L} \right) \sin \left( \frac{\pi y}{L} \right) + \eps U 
\sin \frac{\pi x}{L} \sin \frac{U t}{L}
\end{array}
\right\} \, ,
\label{eq:taylorgreencontrolled}
\end{equation}
\end{small}
in which $ p(y,t) $ is given in (\ref{eq:pyt}).  The restrictions on the parameters are
$ U > 0 $, $ L > 0 $, $ \left| \eps \right| \ll 1 $, $ - \infty < T_s < 0 $ and $ - \infty < S $.
 The $ p $ restriction for describing the manifold 
could alternatively be given
as $ y $-restriction $ 0 < y < y_m(t) $, with $ y_m(t) $ defined in (\ref{eq:taylorgreenymax}). This is only a limitation of $ y $ in describing the $ {\mathcal O}(\eps)-$close
manifold; there is no restriction of $ y $ in the flow (\ref{eq:taylorgreencontrolled}).

\begin{figure}
\begin{tabular}{c}
\includegraphics[width=10cm]{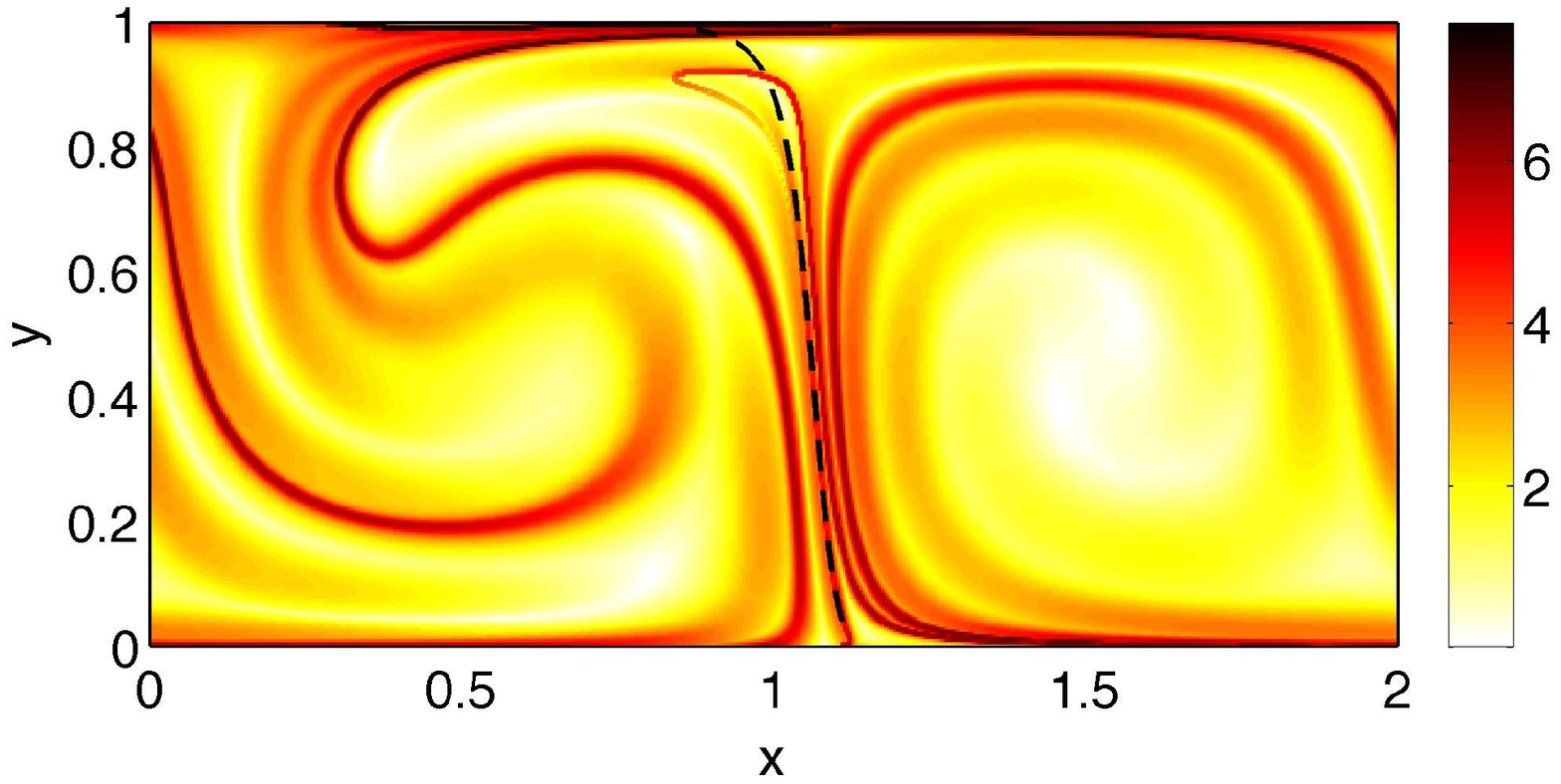}\\
\includegraphics[width=10cm]{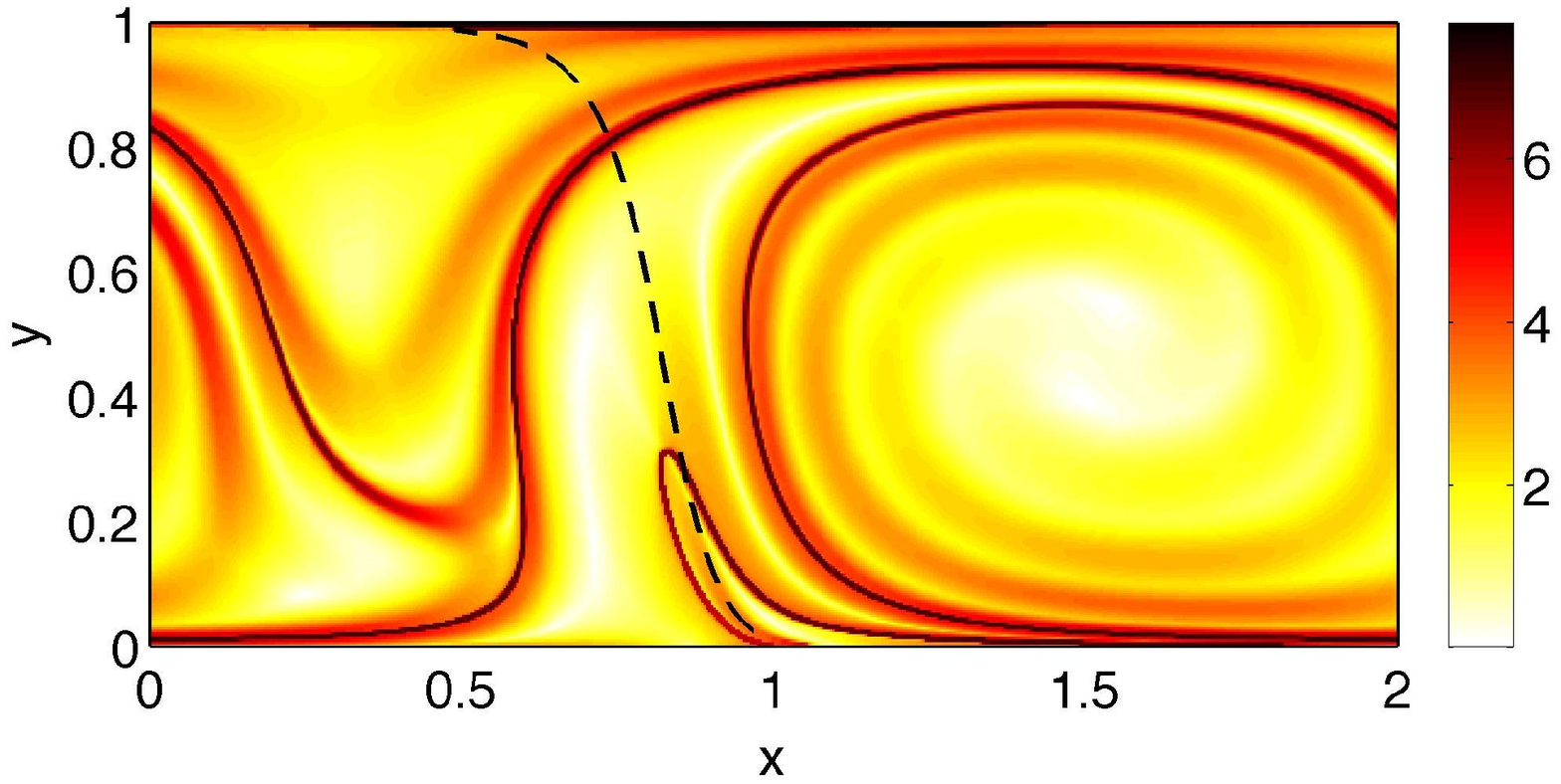}\\
\includegraphics[width=10cm]{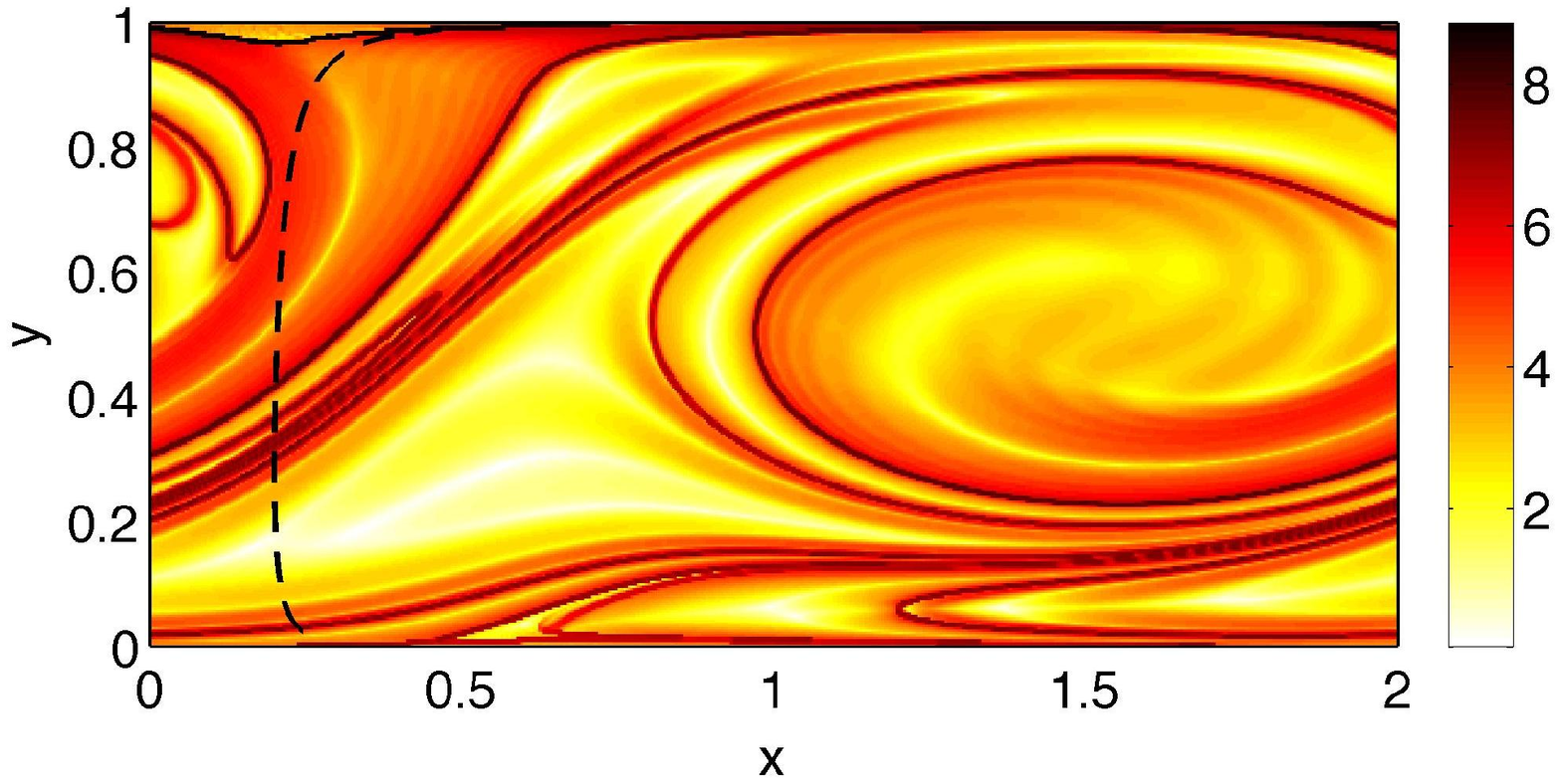}\\
\end{tabular}
\caption{Finite-time Lyapunov exponent fields for (\ref{eq:taylorgreencontrolled}) with $U=1$, $L=1$, $T_s=-1$ and $\eps=0.1$. The desired stable manifold (\ref{eq:taylorgreenmanifold}) is shown by the black dashed curve, and the panels are respectively for the choices $t=0.2, \, 0.5, \, 1.5$.}\label{fig:ftle_Ts}
\end{figure}

In order to test the validity of the analytical results, we compare them with numerically approximated manifolds for the system (\ref{eq:taylorgreencontrolled}). For this we approximate the respective finite-time Lyapunov (FTLE) fields, choosing an integration interval $[t, t+1]$. Ridges in the FTLE field at $t$ indicate---under certain additional assumptions \cite{hallervariational}---the location of stable manifolds. We refer to \cite{padbergplasma} for a brief explanation of the computational scheme used in this paper. Recent work by Haller \cite{hallervariational} sets the heuristical FTLE approach on a sound mathematical basis.

In Fig.\ \ref{fig:ftle_eps} we show the finite-time Lyapunov fields computed for the system (\ref{eq:taylorgreencontrolled}) at $t=-0.9$, with the choice of parameters $U=1$, $L=1$ and $T_s=-1$.   
The black dashed  curve indicates the desired stable manifold (\ref{eq:taylorgreenmanifold}). This desired stable manifold matches up well with a ridge of the FTLE field, in particular, when $\eps$ is sufficiently small. One clearly sees deviations for $y$ being close to $1$ in the two upper panels of Fig.~\ref{fig:ftle_eps}. In the bottom panel, we have chosen a larger value of $\eps$; here the alignment of the desired manifold and the numerically observed one breaks down already for small $y$. The lack of control of the manifold for larger $y$ is a reflection of the condition (\ref{eq:taylorgreenymax}); the $\mathcal{O}(\eps^2)$-closeness of the desired manifold to the true manifold breaks down beyond this value.
 
In contrast, we investigate the worsening of the control strategy with time (at fixed $\eps$) in Fig.~\ref{fig:ftle_Ts}. These and other experiments indicate that the control strategy works well in the range $T_s \le t < T_s + 1.3$ for this example. The reason for the worsening which occurs for larger values of t in this example can be explained by viewing the last panel ($t = 1.5$) in Fig.~\ref{fig:ftle_Ts}. Here, the mappability of the required stable manifold is being compromised near the hyperbolic point along the $y = 0$ line; the black dashed curve is becoming perpendicular to the line $x = 1$ (which is the unperturbed stable manifold). Therefore, the domain $[S_1(t),S_2(t))$ associated with the legitimacy of the control strategy appears to be shrinking at such larger $t$ values.

\section{Concluding remarks}

We have in this article developed a theoretical framework based on which it is possible to 
move a stable/unstable manifold in a two-dimensional autonomous system, to a {\em desired} nonautonomous location which is subjected to certain mappability conditions to 
the original manifold.  A rigorous error estimate for the procedure was developed.  A
numerical example is used to demonstrate the efficacy of the manifold control method.  To our knowledge, this is the first
study which furnishes a method for controlling stable and unstable manifolds nonautonomously in the sense of making
them follow a user-specified time-variation.

In a forthcoming article, we will develop methods for simplifying the hypotheses required
for the restricted stable and unstable manifolds, in order to address the computationally
natural situation of attempting to achieve a desired stable/unstable manifold which is
given in the form $ f(x,y,t) = 0 $, as opposed to having to work through the parameter $ p $.  Preliminary results indicate that
the control strategy can be implemented, for example, to achieve highly wiggly user-specified nonautonomous invariant manifolds.
We expect to obtain insights into a more natural implementation of the mappability condition, so that unreasonable expectations from our control strategy (such as the dashed
curve we tried to require in the final panel in Fig.~\ref{fig:ftle_Ts}) are avoided.
Extensive numerical analyses will be performed in all these situations.

This article complements the authors' work on controlling hyperbolic trajectories
(that is, the ``beginning of stable/unstable manifolds'').  In ongoing research, recent two-dimensional control strategies \cite{control} are being extended to arbitrary dimensions,
and to arbitrarily high-order accuracy.   Building on the present article, similarly extending control strategies to stable/unstable manifolds in high dimensions shall
be our next focus. 

\vspace*{0.4cm}

\noindent
{\bf Acknowledgements:} 
This work was partially supported by a grant from the Simons Foundation (\#236923 to SB),
and by TU Dresden, in sponsoring a visit by SB to Dresden. A start-up grant from the University of
Adelaide to SB is also gratefully acknowledged. \\



\begin{thebibliography}{29}
\expandafter\ifx\csname natexlab\endcsname\relax\def\natexlab#1{#1}\fi
\providecommand{\url}[1]{\texttt{#1}}
\providecommand{\href}[2]{#2}
\providecommand{\path}[1]{#1}
\providecommand{\DOIprefix}{doi:}
\providecommand{\ArXivprefix}{arXiv:}
\providecommand{\URLprefix}{URL: }
\providecommand{\Pubmedprefix}{pmid:}
\providecommand{\doi}[1]{\href{http://dx.doi.org/#1}{\path{#1}}}
\providecommand{\Pubmed}[1]{\href{pmid:#1}{\path{#1}}}
\providecommand{\bibinfo}[2]{#2}
\ifx\xfnm\relax \def\xfnm[#1]{\unskip,\space#1}\fi
\bibitem[{Haller and Beron-Vera(2012)}]{hallerberonvera}
\bibinfo{author}{G.~Haller}, \bibinfo{author}{F.~Beron-Vera},
\newblock \bibinfo{title}{Geodesic theory of transport barriers in
  two-dimensional flows},
\newblock \bibinfo{journal}{Phys. D} \bibinfo{volume}{241}
  (\bibinfo{year}{2012}) \bibinfo{pages}{1680--1702}.
\bibitem[{Blazevski and Haller(2013)}]{blazevskihaller}
\bibinfo{author}{D.~Blazevski}, \bibinfo{author}{G.~Haller},
\newblock \bibinfo{title}{Hyperbolic and elliptic transport barriers in
  three-dimensional unsteady flows}  (\bibinfo{year}{2013}).
  \bibinfo{note}{{a}rXiv: submit/0752347 (submitted)}.
\bibitem[{Allshouse and Thiffeault(2012)}]{allshousethiffeault}
\bibinfo{author}{M.~Allshouse}, \bibinfo{author}{J.-L. Thiffeault},
\newblock \bibinfo{title}{Detecting coherent structures using braids},
\newblock \bibinfo{journal}{Phys. D} \bibinfo{volume}{241}
  (\bibinfo{year}{2012}) \bibinfo{pages}{95--105}.
\bibitem[{Froyland and Padberg(2012)}]{froylandpadbergentropy}
\bibinfo{author}{G.~Froyland}, \bibinfo{author}{K.~Padberg},
\newblock \bibinfo{title}{Finite-time entropy: a probabilistic method for
  measuring nonlinear stretching},
\newblock \bibinfo{journal}{Phys. D} \bibinfo{volume}{241}
  (\bibinfo{year}{2012}) \bibinfo{pages}{1612--1628}.
\bibitem[{Budi\u{s}i\'c and Mezi\'c(2012)}]{budisicmezic}
\bibinfo{author}{M.~Budi\u{s}i\'c}, \bibinfo{author}{I.~Mezi\'c},
\newblock \bibinfo{title}{Geometry of ergodic quotient reveals coherent
  structures in flows},
\newblock \bibinfo{journal}{Phys. D} \bibinfo{volume}{241}
  (\bibinfo{year}{2012}) \bibinfo{pages}{1255--1269}.
\bibitem[{Haller(2011)}]{hallervariational}
\bibinfo{author}{G.~Haller},
\newblock \bibinfo{title}{A variational theory for {L}agrangian {C}oherent
  {S}tructures},
\newblock \bibinfo{journal}{Phys. D} \bibinfo{volume}{240}
  (\bibinfo{year}{2011}) \bibinfo{pages}{574--598}.
\bibitem[{Shadden et~al.(2005)Shadden, Lekien, and Marsden}]{shadden}
\bibinfo{author}{S.~Shadden}, \bibinfo{author}{F.~Lekien},
  \bibinfo{author}{J.~Marsden},
\newblock \bibinfo{title}{Definition and properties of {L}agrangian coherent
  structures from finite-time {L}yapunov exponents in two-dimensional aperiodic
  flows},
\newblock \bibinfo{journal}{Phys. D} \bibinfo{volume}{212}
  (\bibinfo{year}{2005}) \bibinfo{pages}{271--304}.
\bibitem[{Dellnitz and Junge(1997)}]{dellnitzjunge}
\bibinfo{author}{M.~Dellnitz}, \bibinfo{author}{O.~Junge},
\newblock \bibinfo{title}{Almost invariant sets in {C}hua's circuit},
\newblock \bibinfo{journal}{Int. J. Bif. Chaos} \bibinfo{volume}{7}
  (\bibinfo{year}{1997}) \bibinfo{pages}{2475--2485}.
\bibitem[{Froyland et~al.(2010{\natexlab{a}})Froyland, Lloyd, and
  Santitissadeekorn}]{froylandcoherent}
\bibinfo{author}{G.~Froyland}, \bibinfo{author}{S.~Lloyd},
  \bibinfo{author}{N.~Santitissadeekorn},
\newblock \bibinfo{title}{Coherent sets for nonautonomous dynamical systems},
\newblock \bibinfo{journal}{Phys. D} \bibinfo{volume}{239}
  (\bibinfo{year}{2010}{\natexlab{a}}) \bibinfo{pages}{1527--1541}.
\bibitem[{Froyland et~al.(2010{\natexlab{b}})Froyland, Santitissadeekorn, and
  Monahan}]{froylandchaos}
\bibinfo{author}{G.~Froyland}, \bibinfo{author}{N.~Santitissadeekorn},
  \bibinfo{author}{A.~Monahan},
\newblock \bibinfo{title}{Transport in time-dependent dynamical systems:
  finite-time coherent sets},
\newblock \bibinfo{journal}{Chaos} \bibinfo{volume}{20}
  (\bibinfo{year}{2010}{\natexlab{b}}) \bibinfo{pages}{043116}.
\bibitem[{Mezi\'c et~al.(2010)Mezi\'c, Loire, Fonoberov, and
  Hogan}]{mezicscience}
\bibinfo{author}{I.~Mezi\'c}, \bibinfo{author}{S.~Loire},
  \bibinfo{author}{V.~Fonoberov}, \bibinfo{author}{P.~Hogan},
\newblock \bibinfo{title}{A new mixing diagnostic and {G}ulf oil spill
  movement},
\newblock \bibinfo{journal}{Science} \bibinfo{volume}{330}
  (\bibinfo{year}{2010}) \bibinfo{pages}{486--489}.
\bibitem[{Balasuriya(2011)}]{tangential}
\bibinfo{author}{S.~Balasuriya},
\newblock \bibinfo{title}{A tangential displacement theory for locating
  perturbed saddles and their manifolds},
\newblock \bibinfo{journal}{SIAM J. Appl. Dyn. Sys.} \bibinfo{volume}{10}
  (\bibinfo{year}{2011}) \bibinfo{pages}{1100--1126}.
\bibitem[{Balasuriya(2013)}]{open}
\bibinfo{author}{S.~Balasuriya},
\newblock \bibinfo{title}{Nonautonomous flows as open dynamical sytems:
  characterising escape rates and time-varying boundaries},
\newblock in: \bibinfo{booktitle}{Ergodic Theory, Open Dynamics and
  Structures}, \bibinfo{publisher}{Springer}, \bibinfo{year}{2013}.
  \bibinfo{note}{In press}.
\bibitem[{Balasuriya(2012)}]{unsteady}
\bibinfo{author}{S.~Balasuriya},
\newblock \bibinfo{title}{Explicit invariant manifolds and specialised
  trajectories in a class of unsteady flows},
\newblock \bibinfo{journal}{Phys. Fluids} \bibinfo{volume}{24}
  (\bibinfo{year}{2012}) \bibinfo{pages}{127101}.
\bibitem[{Farazmand and Haller(2012)}]{farazmandhaller}
\bibinfo{author}{M.~Farazmand}, \bibinfo{author}{G.~Haller},
\newblock \bibinfo{title}{Computing {L}agrangian {C}oherent {S}tructures from
  variational {LCS} theory},
\newblock \bibinfo{journal}{Chaos} \bibinfo{volume}{22} (\bibinfo{year}{2012})
  \bibinfo{pages}{013128}.
\bibitem[{Peacock and Dabiri(2010)}]{peacockdabiri}
\bibinfo{author}{T.~Peacock}, \bibinfo{author}{J.~Dabiri},
\newblock \bibinfo{title}{Introduction to focus issue: {L}agrangian {C}oherent
  {S}tructures},
\newblock \bibinfo{journal}{Chaos} \bibinfo{volume}{20} (\bibinfo{year}{2010})
  \bibinfo{pages}{017501}.
\bibitem[{Boffetta et~al.(2001)Boffetta, Lacorata, Radaelli, and
  Vulpiani}]{boffetta}
\bibinfo{author}{G.~Boffetta}, \bibinfo{author}{G.~Lacorata},
  \bibinfo{author}{G.~Radaelli}, \bibinfo{author}{A.~Vulpiani},
\newblock \bibinfo{title}{Detecting barriers to transport: a review of
  different techniques},
\newblock \bibinfo{journal}{Phys. D} \bibinfo{volume}{159}
  (\bibinfo{year}{2001}) \bibinfo{pages}{58--70}.
\bibitem[{Krauskopf et~al.(2005)Krauskopf, Osinga, Doedel, Henderson,
  Guckenheimer, Vladimirsky, Dellnitz, and Junge}]{krauskopf}
\bibinfo{author}{B.~Krauskopf}, \bibinfo{author}{H.~Osinga},
  \bibinfo{author}{E.~Doedel}, \bibinfo{author}{M.~Henderson},
  \bibinfo{author}{J.~Guckenheimer}, \bibinfo{author}{A.~Vladimirsky},
  \bibinfo{author}{M.~Dellnitz}, \bibinfo{author}{O.~Junge},
\newblock \bibinfo{title}{A survey of method's for computing (un)stable
  manifold of vector fields},
\newblock \bibinfo{journal}{Int. J. Bif. Chaos} \bibinfo{volume}{15}
  (\bibinfo{year}{2005}) \bibinfo{pages}{763--791}.
\bibitem[{Froyland(2013)}]{froylandanalytic}
\bibinfo{author}{G.~Froyland},
\newblock \bibinfo{title}{An analytic framework for identifying finite-time
  coherent sets in time-dependent dynamical systems},
\newblock \bibinfo{journal}{Phys. D} \bibinfo{volume}{250}
  (\bibinfo{year}{2013}) \bibinfo{pages}{1--19}.
\bibitem[{Radko(2011)}]{radko}
\bibinfo{author}{T.~Radko},
\newblock \bibinfo{title}{On the generation of large-scale structures in a
  homogeneous eddy field},
\newblock \bibinfo{journal}{J. Fluid Mech.} \bibinfo{volume}{668}
  (\bibinfo{year}{2011}) \bibinfo{pages}{76--99}.
\bibitem[{Chandrasesekhar(1961)}]{chandrasekhar}
\bibinfo{author}{S.~Chandrasesekhar}, \bibinfo{title}{Hydrodynamics and
  hydrodynamic stability}, \bibinfo{publisher}{Dover}, \bibinfo{address}{New
  York}, \bibinfo{year}{1961}.
\bibitem[{Balasuriya(2005{\natexlab{a}})}]{mixer}
\bibinfo{author}{S.~Balasuriya},
\newblock \bibinfo{title}{Approach for maximizing chaotic mixing in
  microfluidic devices},
\newblock \bibinfo{journal}{Phys. Fluids} \bibinfo{volume}{17}
  (\bibinfo{year}{2005}{\natexlab{a}}) \bibinfo{pages}{118103}.
\bibitem[{Balasuriya(2005{\natexlab{b}})}]{periodic}
\bibinfo{author}{S.~Balasuriya},
\newblock \bibinfo{title}{Direct chaotic flux quantification in perturbed
  planar flows: general time-periodicity},
\newblock \bibinfo{journal}{SIAM J. Appl. Dyn. Sys.} \bibinfo{volume}{4}
  (\bibinfo{year}{2005}{\natexlab{b}}) \bibinfo{pages}{282--311}.
\bibitem[{Balasuriya and Finn(2012)}]{l2mixer}
\bibinfo{author}{S.~Balasuriya}, \bibinfo{author}{M.~Finn},
\newblock \bibinfo{title}{Energy constrained transport maximization across a
  fluid interface},
\newblock \bibinfo{journal}{Phys. Rev. Lett.} \bibinfo{volume}{108}
  (\bibinfo{year}{2012}) \bibinfo{pages}{244503}.
\bibitem[{Coppel(1978)}]{coppel}
\bibinfo{author}{W.~A. Coppel}, \bibinfo{title}{Dichotomies in Stability
  Theory}, number \bibinfo{number}{629} in \bibinfo{series}{Lecture Notes in
  Mathematics}, \bibinfo{publisher}{Springer-Verlag},
  \bibinfo{address}{Berlin}, \bibinfo{year}{1978}.
\bibitem[{Yi(1993)}]{yi}
\bibinfo{author}{Y.~Yi},
\newblock \bibinfo{title}{A generalized integral manifold theorem},
\newblock \bibinfo{journal}{J. Differential Equations} \bibinfo{volume}{102}
  (\bibinfo{year}{1993}) \bibinfo{pages}{153--187}.
\bibitem[{Yagasaki(2008)}]{yagasakicontrol}
\bibinfo{author}{K.~Yagasaki},
\newblock \bibinfo{title}{Invariant manifolds and control of hyperbolic
  trajectories on infinite- or finite-time intervals},
\newblock \bibinfo{journal}{Dyn. Sys.} \bibinfo{volume}{23}
  (\bibinfo{year}{2008}) \bibinfo{pages}{309--331}.
\bibitem[{Padberg et~al.(2007)Padberg, Hauff, Jenko, and Junge}]{padbergplasma}
\bibinfo{author}{K.~Padberg}, \bibinfo{author}{T.~Hauff},
  \bibinfo{author}{F.~Jenko}, \bibinfo{author}{O.~Junge},
\newblock \bibinfo{title}{Lagrangian structures and transport in turbulent
  magnetized plasmas},
\newblock \bibinfo{journal}{New J. Phys.} \bibinfo{volume}{9}
  (\bibinfo{year}{2007}) \bibinfo{pages}{400}.
\bibitem[{Balasuriya and Padberg-{G}ehle(2013)}]{control}
\bibinfo{author}{S.~Balasuriya}, \bibinfo{author}{K.~Padberg-{G}ehle},
\newblock \bibinfo{title}{Controlling the unsteady analogue of saddle
  stagnation points},
\newblock \bibinfo{journal}{SIAM J.\ Appl.\ Math.} \bibinfo{volume}{73}
  (\bibinfo{year}{2013}) \bibinfo{pages}{1038--1057}.

\end{thebibliography}
\end{document}